\documentclass[11pt,reqno]{amsart}

\title[]{On the Global Well-Posedness and Gevrey Regularity of some Electrodiffusion Models}

\author{Elie Abdo}
\address{Department of Mathematics, Temple University, Philadelphia, PA 19122, USA}
\email{elie.abdo@temple.edu}
\author{Fizay-Noah Lee}
\address{Program in applied and computational mathematics, Princeton University, Princeton, NJ 08544, USA}
\email{fl6@math.princeton.edu}
\author{Weinan Wang}
\address{Department of Mathematics, University of Arizona, Tucson, AZ, 85712, USA}
\email{weinanwang@math.arizona.edu}

\usepackage[margin=1in]{geometry}
\usepackage{amsmath, amsthm, amssymb}
\usepackage{times}
\usepackage{color}
\usepackage{hyperref}
\newcommand{\pa}{\partial}
\newcommand{\la}{\label}
\newcommand{\fr}{\frac}
\newcommand{\na}{\nabla}
\newcommand{\be}{\begin{equation}}
\newcommand{\ee}{\end{equation}}
\newcommand{\ba}{\begin{array}{l}}
	\newcommand{\ea}{\end{array}}

\newcommand{\beg}{\begin}

\newcommand{\D}{\Delta}
\renewcommand{\l}{\Lambda}

\usepackage{MnSymbol,wasysym}

\newcommand{\R}{\mathbb R}

\def\ZZ{{\mathbb Z}}
\def\RR{{\mathbb R}}
\def\TT{{\mathbb T}}

\def\PP{\mathbb P}

\date{\today}
\begin{document}
	\begin{abstract} We consider the Nernst-Planck equations describing the nonlinear time evolution of multiple ionic concentrations in a two-dimensional incompressible fluid.  The velocity of the fluid evolves according to either the Euler or Darcy's  equations, both forced nonlinearly by the electric forces generated by the presence of charged ions. We address the global well-posedness and Gevrey regularity of the resulting electrodiffusion models in the periodic setting.
	\end{abstract}

	\maketitle
	\tableofcontents
	\section{Introduction}\la{intro}

		Electrodiffusion of ions in a fluid is governed by three main mechanisms: diffusion driven by the gradient of the ionic concentrations, transport driven by the gradient of the electrical potential due to the presence of ions, and transport due to the fluid. Mathematically, electrodiffusion is modeled by the Nernst-Planck equations, which are mass-balance equations describing the time evolution of local ionic concentrations. The equations relate the ionic fluxes to diffusion and transport by the electrical field and the fluid  velocity. In turn, the electrical field is determined nonlocally via a Poisson equation from the ionic concenrations. Real world applications motivate studies of electrodiffusion phenomena. One such model is addressed in \cite{TR} in the context of improving the performance and lifetime of batteries, which has direct applications to the development of electric vehicles, portable electronics, etc. Other applications of electrodiffusion systems arise in semiconductors \cite{gaj,mock} and ion selective membranes \cite{davidson}.

		The \textit{Nernst-Planck equations} are given by
	\be \la{1}
	\pa_t c_i + u \cdot \na c_i - D_i \Delta c_i = D_i z_i \na \cdot (c_i \na \Phi), \quad i=1,...,n
	\ee describing the time evolution of the ionic concentrations $c_i$  of $n$ ionic species, with valences $z_i\in\mathbb{R}$ and diffusivities $D_i>0$ in a two-dimensional incompressible fluid. Above, the electrical potential $\Phi$ is determined by the Poisson equation
	\be \la{4}
	- \Delta \Phi = \rho = \sum\limits_{i=1}^{n} z_ic_i.
	\ee
	
	In this paper, we consider two electrodiffusion models distinguished by the choice of model used to describe the fluid velocity: the velocity $u$ and pressure $p$ of the fluid obey either the \textit{Euler equations} 
	\be \la{2}
	\pa_t u + u \cdot \na u + \na p = - \rho \na \Phi,\quad \na\cdot u=0
	\ee or \textit{Darcy's law}
	\be \la{3}
	u + \na p = - \rho \na \Phi,\quad \na\cdot u=0.
	\ee 
	
	The model described by \eqref{1}, \eqref{4} and \eqref{2} is called the \textit{Nernst-Plank-Euler system} and will be referred to as NPE, whereas the model described by \eqref{1}, \eqref{4} and \eqref{3} is called the\textit{ Nernst-Planck-Darcy system} and will be referred to as NPD. Physically, considering Euler equations amounts to considering a region of low viscosity fluid far from any solid boundaries. On the other hand, Darcy's law corresponds to creeping flow in a porous medium, which is of physical interest, for example, when considering electrodiffusion in porous electrodes. In this paper, we study the global well-posedness and Gevrey regularity of the NPE and NPD systems in the periodic setting on the two-dimensional torus $\TT^2 = [0, 2\pi]^2$. 
	
	The global well-posedness of the periodic NPE and NPD systems were addressed respectively in \cite{IS} and \cite{IS1} for two ionic species having valences 1 and -1 and equal diffusivities. The special structure of the two-species NPE and NPD models yields energy estimates that facilitate the analysis. These same bounds also yield precise information on the long time behavior of solutions. In fact, it is shown in \cite{IS1} and \cite{IS} that the the two ionic concentrations converge in $L^2$ towards constant distributions (corresponding to their initial spatial averages) exponentially quickly in time. In this paper, we are interested in the more general models that model the time evolution of $n$ ionic species with different valences and diffusivities. At this level of generality, the ideas of \cite{IS1} and \cite{IS} no longer apply. Instead, in this paper we adapt ideas from \cite{AI6} and \cite{CI} to obtain time dependent $L^2$ bounds for the velocity and ionic concentrations, and we use them to bootstrap to obtain higher regularity bounds, assuming sufficiently regular initial data. We thus establish global well-posedness of solutions of NPE and NPD in all sufficiently regular Sobolev spaces.
	
	We also consider the analyticity of solutions to the NPE and NPD models. Gevrey regularity techniques were initiated by Foias and Temam in \cite{FT} to study the analyticity of solutions of the Navier-Stokes equations. Their idea uses Fourier series expansions, is restricted to Hilbert spaces, and requires Sobolev regular initial regularity. In \cite{GK}, Grujic and Kukavica introduced a simpler method to obtain the analyticity of the Navier-Stokes system for $L^p$ initial data based on a mild formulation of the complexified problem.  Bardos and Benachour proved in \cite{BB} the persistence of analyticity for the Euler equations on three-dimensional bounded domains and obtained decaying rates for the radius of analyticity. Using the method of Gevrey regularity, Levermore and Oliver obtained the same rates in \cite{LO} for the three-dimensional periodic Euler equations. In \cite{KV}, Kukavica and Vicol improved the aforementioned results and derived lower bounds for the radius of analyticity for solutions of the three-dimensional Euler system that depend algebraically on $\exp\int_0^t\|\na u(s)\|_{L^\infty}\,ds$. 
	
	The NPE and NPD models are closely related to the Nernst-Planck-Navier-Stokes (NPNS) system, where the velocity of the fluid evolves according to Navier-Stokes equations forced nonlinearly by the electrical force $-\rho \na \Phi$. The NPNS system was investigated on bounded domains in \cite{CI}, \cite{cil1}, \cite{cil2}, \cite{cil3}, \cite{cil4}, and on the 2D periodic torus in \cite{AI4}. In \cite{AI6}, the authors adapted the technique of \cite{GK} to study the analyticity of the periodic NPNS system, forced by body forces in the fluid and additional charge densities, for $L^p$ initial data. However, the idea of the proof applies to semi-linear parabolic equations and breaks down in the absence of the fluid dissipation. In this paper, we use ideas from \cite{FT} and \cite{KV} and provide quantitative estimates for the growth in time of the Gevrey norms of the solution to both the NPE and NPD systems. These estimates remain finite on all finite intervals, and thus our results imply global in time spatial analyticity of solutions.

	The Nernst-Planck equations have a dissipative structure, which is taken advantage of to estimate the nonlinearities involving ionic concentration terms. The Euler equations, on the other hand, are inviscid and thus nondissipative. Despite the inviscidity of the Euler system and the nonlinearties driving the time evolution of the concentrations and velocity, we successfully adapt the ideas of \cite{KV} and establish the analyticity of solutions of NPE via use of Fourier series techniques. The challenges arising from the analysis of the NPD system are different though. NPD is effectively a system of $n$ dissipative partial differential equations with $n$ unknowns, namely the ionic concentrations $c_1, ..., c_n$, and the fluid velocity is computed from the charge density $\rho$ via the nonlinear relation 
	\be 
	u = - \PP (\rho \na \Phi)
	\ee where $\PP$ is the Leray-Hodge projector onto the space of divergence-free vector fields. Hence, the nonlinear transport terms $u \cdot \na c_i$ are effectively cubic. Controlling these terms is the main source of difficulty and is dealt with using Fourier series techniques, yielding bounds with cubic dependence on the desired Gevrey norm. This gives local analytic solutions, which we then extend globally by establishing uniform control of high regularity Sobolev norms of the solution.

	This paper is organized as follows. In Section \ref{pre}, we introduce the functional setting and notation conventions that will be frequently used throughout the paper. In Section \ref{S2}, we study the global well-posedness of the NPE system on $\TT^2 \times [0,T]$. Indeed, we first prove the existence of a unique local regular solution provided that the initial velocity and ionic concentrations are at least Sobolev $H^3$ regular, and then we show that this local solution extends to a global regular solution on the time interval $[0,T]$. In Section \ref{S3}, we show that the solution to the NPE system is Gevrey regular with a radius of analyticity that depends on the Sobolev regularity of the solution on the whole time interval $[0,T]$ under the assumption that the initial data is real analytic. In Section \ref{S4}, we show the global well-posedness of the NPD system on $\TT^2 \times [0,T]$ and we prove the existence of a global unique regular solution for at least Sobolev $H^{\fr{3}{2}}$ regular initial concentrations. Finally, we show in Section \ref{S5} that the solution to the NPD system is space analytic, for any initial data that is more regular that Sobolev $H^2$, with a radius of analyticity depending only on the initial regularity of the velocity and concentrations. 
	
\section{Functional Settings and Notations} \la{pre}
	
For $1 \le p \le \infty$, we denote by $L^p(\TT^2)$ the Lebesgue spaces of measurable periodic functions $f$ from $\TT^2$ to $\R$ (or $\RR^2)$ such that 
\be 
\|f\|_{L^p} = \left(\int_{\TT^2} \|f\|^p dx\right)^{1/p} <\infty
\ee if $p \in [1, \infty)$ and 
\be 
\|f\|_{L^{\infty}} = {\mathrm{esssup}}_{\TT^2}  |f| < \infty
\ee if $p = \infty$. The $L^2(\TT^2)$ inner product is denoted by $(\cdot,\cdot)_{L^2}$.

For $s\in\R$, the fractional Laplacian $\l^s$ applied to a mean zero function $f$ is defined as a Fourier multiplier with symbol $|k|^s$, that is, for $f$ given by
\be
f = \sum\limits_{k \in \mathbb{Z}^2 \setminus \left\{0\right\}} f_k e^{ik \cdot x},
\ee and obeying 
\be 
\sum\limits_{k \in \mathbb{Z}^2 \setminus \left\{0\right\}} |k|^{2s} |f_k|^2 < \infty,
\ee
we define
\be
\l^s f = \sum\limits_{k \in \mathbb{Z}^2 \setminus \left\{0\right\}} |k|^s f_k e^{ik \cdot x}.
\ee 
For $\tau> 0$, $s > 0$, and functions $f$ obeying 
\be 
\sum\limits_{k \in \mathbb{Z}^2 \setminus \left\{0\right\}} e^{2\tau |k|^{s}} |f_k|^2 < \infty,
\ee we define 
\be 
e^{\tau \l^s} f = \sum\limits_{k \in \mathbb{Z}^2 \setminus \left\{0\right\}} e^{\tau |k|^s} f_k e^{ik \cdot x}.
\ee

We denote by $\mathbb{P}$ the Leray-Hodge projection onto the space divergence free vector fields. For a mean-free periodic vector field $v = (v_1, v_2)$ with Fourier series
\be 
v = \sum\limits_{j \in \mathbb{Z}^2 \setminus \left\{0\right\}} v_j e^{ij \cdot x},
\ee  $\PP v$ has the following Fourier representation 
\be
\PP v = \sum_{j\in\ZZ^2\setminus\{0\}} \left[v_j - (v_j\cdot j)\frac{j}{|j|^2} \right] e^{ij\cdot x}.
\la{pu}
\ee The operator $\PP$ is bounded on $L^p$ spaces for any $p \in (1, \infty)$.

For $s > 0$, we denote by $H^s(\TT^2)$ the Sobolev spaces of measurable periodic functions $f$ from $\TT^2$ to $\R$ (or $\RR^2)$ obeying 
\be 
\|f\|_{H^s}^2 = \sum\limits_{k \in \ZZ^2} (1 + |k|^s)^2|f_k|^2 < \infty.
\ee

For a Banach space $(X, \|\cdot\|_{X})$ and $p\in [1,\infty]$, we consider the Lebesgue spaces $ L^p(0,T; X)$ of functions $f$ from $X$ to $\R$ (or $\RR^2)$ satisfying 
\be 
\int_{0}^{T} \|f\|_{X}^p dt  <\infty
\ee with the usual convention when $p = \infty$. 

Throughout the paper, $C$ denotes a positive constant depending only on the parameters of the problem (namely the diffusivities and valences of the ionic concentrations) and some universal constants, and it changes from line to line along the proofs. The dependency on any other variable will be stated explicitly. We also use the notation convention $A\lesssim B$ when $A \leq CB$. Finally, we adopt the notation $[\Lambda^m, v\cdot\na]w$ to denote the commutator $\l^{m} (v \cdot \na w) - v \cdot \na \l^{m}w$.
	
\section{Global Well-posedness of the Nernst-Planck-Euler system }\la{S2}

In this section, we address the global existence and uniqueness of solutions to  the NPE system:

\beg{thm} \la{T1} (Global well-posedness and regularity of NPE) Let $T>0$ be arbitrary and $m \ge 3$.  Assume that  the initial ionic concentrations $c_i(0) \in H^m$ are nonnegative, and the initial velocity $u_0 \in H^m$ is mean-zero and divergence-free. Then the NPE system  described by \eqref{1}, \eqref{4} and \eqref{2} has a unique solution $(u, c_1, ..., c_n)$ on the time interval $[0,T]$ with the property that
\be 
(u, c_1, ..., c_n) \in (L^{\infty}(0,T; H^m))^{n+1}.
\ee
\end{thm}

The proof of Theorem \ref{T1} will be presented in this section. First, we note that the spatial averages of the ionic concentrations are constant in time, that is 
\be 
\int c_i(t) dx = \int c_i(0) dx
\ee for all $t \ge 0$ and all $i \in \left\{1, ..., n\right\}$. Moreover, the spatial average of the velocity vanishes at any positive time $t$, a fact that follows from integrating the forced Euler equation \eqref{2} obeyed by $u$ and noting that 
\be
\beg{aligned} 
&\int \rho \na \Phi dx 
= \int \rho \na \l^{-2} \rho dx 
= \int \l^{-2} \rho \na  \rho dx 
= - \int (\na \l^{-2} \rho) \rho dx
= - \int \rho \na \Phi dx
\end{aligned}
\ee from which we conclude that the spatial integral of $\rho \na \Phi$ vanishes at any time $t \ge 0$. 

Now we proceed to prove Theorem \ref{T1}. The proof is divided into several steps.

\textbf{Step 1. Local well-posedness for $H^m$ initial datum.} We sketch a proof of local existence, with a focus on obtaining a priori estimates. Rigorously, we may do the following computations for a mollified version of the  NPE system \cite{MB}. We fix $m>2$, and consider initial conditions in $H^m$. Then, we apply $\Lambda^m$ to \eqref{1}, multiply the resulting equation by $\Lambda^m c_i$, and integrate by parts:
\be
\begin{aligned}\la{le1}
\fr{1}{2}\fr{d}{dt}\|\Lambda^m c_i\|_{L^2}^2+D_i\|\Lambda^m\na c_i\|_{L^2}^2=-\left(\Lambda^m(u\cdot\na c_i),\Lambda^m c_i\right)_{L^2}+D_iz_i\left(\Lambda^m\na\cdot(c_i\na\Phi),\Lambda^m c_i\right)_{L^2}.
\end{aligned}
\ee
Due to the fact that $\na\cdot u=0$, we have that $(\Lambda^m(u\cdot\na c_i),\Lambda^mc_i)_{L^2}=([\Lambda^m,u\cdot\na]c_i,\Lambda^m c_i)_{L^2}$. So it follows from commutator estimates \cite{KPV} that
\be
\begin{aligned}\la{le2}
\left|\left([\Lambda^m,u\cdot\na]c_i,\Lambda^m c_i\right)_{L^2}\right|\lesssim& (\|\Lambda^m u\|_{L^2}\|\na c_i\|_{L^\infty}+\|\na u\|_{L^\infty}\|\Lambda^m c_i\|_{L^2})\|\Lambda^m c_i\|_{L^2}\\
\lesssim& \|\Lambda^m u\|_{L^2}\|\Lambda^m c_i\|_{L^2}^2\\
\lesssim& 1+\|\Lambda^m u\|_{L^2}^4+\|\Lambda^m c_i\|_{L^2}^4,
\end{aligned}
\ee
where the second line follows from the embedding $H^m\hookrightarrow W^{1,\infty}$ which holds for $m>2$. We estimate the second term on the right hand side of \eqref{le1} as follows
\be
\begin{aligned}\la{le3}
\left|\left(\Lambda^m\na\cdot(c_i\na\Phi),\Lambda^m c_i\right)_{L^2}\right|=&\left|\left(\Lambda^m(c_i\na\Phi),\Lambda^m \na c_i\right)_{L^2}\right|\\
\lesssim&\|\Lambda^m c_i\|_{L^2}\|\Lambda^m \na \Phi\|_{L^2}\|\Lambda^m \na c_i\|_{L^2} \\
\le& \fr{D_i}{2}\|\Lambda^m\na c_i\|_{L^2}^2+C\sum_{j=1}^n\|\Lambda^m c_j\|_{L^2}^4.
\end{aligned}
\ee
Above, we used the fact that for $m>1$, $H^m$ is a Banach algebra so that $\|\l^m(c_i\na\Phi)\|_{L^2}\lesssim \|\l^m c_i\|_{L^2}\|\l^m\na\Phi\|_{L^2}$. Putting together \eqref{le1}-\eqref{le3}, and summing in $i$, we have thus far
\be\la{lmc}
\fr{d}{dt} \sum\limits_{i=1}^{n}\|\Lambda^m c_i\|_{L^2}^2+ \sum\limits_{i=1}^{n} D_i\|\Lambda^m\na c_i\|_{L^2}^2\lesssim 1+ \|\Lambda^m u\|_{L^2}^4+  \sum\limits_{i=1}^{n} \|\Lambda^m c_i\|_{L^2}^4.
\ee
Next, we  apply $\Lambda^m$ to \eqref{2}, multiply by $\Lambda^m u$, and integrate by parts to obtain
\be\la{le4}
\fr{1}{2}\fr{d}{dt}\|\Lambda^m u\|_{L^2}^2=-\left(\Lambda^m(u\cdot \na u),\Lambda^m u\right)_{L^2}-\left(\Lambda^m(\rho\na\Phi),\Lambda^m u\right)_{L^2}.
\ee
We estimate the terms on the right hand side of \eqref{le4}
\be 
\begin{aligned}
\left|\left(\Lambda^m(u\cdot \na u),\Lambda^m u\right)_{L^2}\right|=&\left|\left([\Lambda^m,u\cdot\na]u,\Lambda^mu\right)_{L^2}\right|
\\ \lesssim&\left(\|\Lambda^m u\|_{L^2}\|\na u\|_{L^\infty}+ \|\na u\|_{L^\infty}\|\Lambda^mu\|_{L^2}\right)\|\Lambda^m u\|_{L^2}
\\ \lesssim& 1+\|\Lambda^m u\|_{L^2}^4       \la{le5}
\end{aligned}
\ee  and 
\be 
\begin{aligned}
\left|\left(\Lambda^m(\rho\na\Phi),\Lambda^m u\right)_{L^2}\right|\lesssim& \|\Lambda^m\rho\|_{L^2}\|\Lambda^m\na\Phi\|_{L^2}\|\Lambda^m u\|_{L^2}\\
\lesssim& 1+\|\Lambda^m u\|_{L^2}^4+  \sum\limits_{i=1}^{n} \|\Lambda^m c_i\|_{L^2}^4.     \la{le6}
\end{aligned}
\ee
Again, in obtaining \eqref{le5}, we used the embedding $H^m\hookrightarrow W^{1,\infty}$ for $m>2$. Putting together \eqref{le4}-\eqref{le6} with \eqref{lmc}, we finally have that 
\be
\fr{d}{dt}\left( \sum\limits_{i=1}^{n} \|\Lambda^m c_i\|_{L^2}^2+\|\Lambda^m u\|_{L^2}^2\right)\lesssim 1+\left( \sum\limits_{i=1}^{n} \|\Lambda^m c_i\|_{L^2}^2+\|\Lambda^m u\|_{L^2}^2\right)^2.
\ee
Then from a Gronwall inequality, this last inequality gives us a time $T^*$, depending on the $H^m$ norms of the initial data, and a local solution of  NPE on $[0,T^*)$ such that
\be\la{d}
\sup_{t\in[0,T^*)}\left( \sum\limits_{i=1}^{n} \|\l^m c_i(t)\|_{L^2}^2+\|\l^m u(t)\|_{L^2}^2\right)\le 2\left( \sum\limits_{i=1}^{n} \|\l^m c_i(0)\|_{L^2}^2+\|\l^m u(0)\|_{L^2}^2\right).
\ee

Uniqueness follows from similar energy estimates. Suppose $(u, c_1,..., c_n)$ and $(u, \bar c_1,...,\bar c_n)$ are two local $H^m$ solutions to  NPE on some common interval $[0,T^*)$. Without loss of generality, we assume that $T^*$ is taken small enough so that both solutions satisfy the local doubling inequality \eqref{d}. Then, the differences $\tilde c_i=c_i-\bar c_i$, $\tilde\Phi=\Phi-\bar\Phi$, $\tilde u=u-\bar u$ satisfy
\begin{align}
    \pa_t \tilde c_i+u\cdot \na\tilde c_i-D_i\D\tilde c_i=&D_iz_i\na\cdot(c_i\na\tilde\Phi+\tilde c_i\na\bar\Phi)-\tilde u\cdot\na\bar c_i \la{npdiff}\\
    \pa_t\tilde u +u\cdot\na\tilde u-\na\tilde p=&-\rho\na\tilde\Phi-\tilde\rho\na\bar\Phi-\tilde u\cdot\na \bar u.    \la{euleridiff}
\end{align}
We multiply the above equations by $\tilde c_i$ and $\tilde u$ respectively, and integrate by parts to obtain
\be
\begin{aligned}\la{tc}
    \fr{1}{2}\fr{d}{dt}\|\tilde c_i\|_{L^2}^2+D_i\|\na\tilde c_i\|_{L^2}^2\lesssim&(\|c_i\|_{L^\infty}\|\na\tilde\Phi\|_{L^2}+\|\tilde c_i\|_{L^2}\|\na\bar\Phi\|_{L^\infty}+\|\tilde u\|_{L^2}\|\bar c_i\|_{L^\infty})\|\na \tilde c_i\|_{L^2}\\
    \le&\fr{D_i}{2}\|\na\tilde c_i\|_{L^2}^2+G_1(t)\left(\sum_{j=1}^n \|\tilde c_j\|_{L^2}^2+\|\tilde u\|_{L^2}^2\right)
\end{aligned}
\ee
and
\be
\begin{aligned}\la{tu}
\fr{1}{2}\fr{d}{dt}\|\tilde u\|_{L^2}^2\lesssim& (\|\rho\|_{L^\infty}\|\na\tilde\Phi\|_{L^2}+\|\tilde\rho\|_{L^2}\|\na\bar\Phi\|_{L^\infty}+\|\tilde u\|_{L^2}\|\na\bar u\|_{L^\infty})\|\tilde u\|_{L^2}\\
\lesssim&G_2(t)\left(\sum_{j=1}^n \|\tilde c_j\|_{L^2}^2+\|\tilde u\|_{L^2}^2\right)
\end{aligned}
\ee
where 
\begin{align}
    G_1(t)=&C\left(\sum_{j=1}^n(\|c_j\|_{L^\infty}^2+\|\bar c_j\|_{L^\infty}^2)+\|\na\bar\Phi\|_{L^\infty}^2\right)\la{G1}\\
    G_2(t)=&\|\rho\|_{L^\infty}+\|\na\bar\Phi\|_{L^\infty}+\|\na\bar u\|_{L^\infty}.
\end{align}
Both $G_1(t)$ and $G_2(t)$ are controlled by the $H^m$ norms of $c_i, \bar c_i,u, \bar u$ and thus integrable on $[0,T^*)$ due to \eqref{d}. Thus, adding \eqref{tc} to \eqref{tu} and integrating, we find
\be
\sup_{t\in[0,T^*)}\left( \sum\limits_{i=1}^{n} \|\tilde c_i(t)\|_{L^2}^2+\|\tilde u(t)\|_{L^2}^2\right)\lesssim \left(\sum\limits_{i=1}^{n} \|\tilde c_i(0)\|_{L^2}^2+\|\tilde u(0)\|_{L^2}^2\right)e^{C\int_0^{T^*} G_1(s)+G_2(s)\,ds}.
\ee
The latter inequality proves uniqueness and continuous (in $L^2$) dependence on initial conditions.

\textbf{Step 2. Positivity of the ionic concentrations.} Suppose $(u, c_1, ..., c_n)$ is a solution of the NPE problem with the property that each ionic concentration $c_i$ belongs to the Lebesgue space $L^2(0,T; H^2)$. Then the positivity of $c_i(0)$ for $i \in \left\{1, ..., n\right\}$ is preserved for all positive times, that is $c_i(x,t) \ge 0$ for all $i \in \left\{1, ..., n\right\}$, for a.e. $x \in \TT^2$, and for all $t \in [0,T]$. The proof can be found in \cite{CI}. 

\textbf{Step 3. A priori $L^2$ uniform-in-time bounds.} Suppose $(u, c_1, ..., c_n)$ is a solution of the NPE problem on the time interval $[0,T]$ such that $c_i(x,t) \ge 0$ for all $i \in \left\{1, ..., n\right\}$, for a.e. $x \in \TT^2$, and for all $t \in [0,T]$. Then there is a positive constant $\Gamma$ depending on the initial data, the time $T$, the parameters of the problem, and some universal constants, such that the following bound
\be \la{Gamma}
\sup\limits_{0 \le t \le T} \left(\|u(t)\|_{L^2} + \sum\limits_{i=1}^{n} \|c_i(t)\|_{L^2} \right) + \int_{0}^{T} \sum\limits_{i=1}^{n} D_i \|\na c_i(t)\|_{L^2}^2 dt \le \Gamma
\ee holds. This a priori bound is proved in Proposition 2 of \cite{AI6}. We omit the details here. 

\textbf{Step 4. A priori $H^1$ uniform-in-time bounds for the velocity $u$.}
The evolution of the vorticity $\omega = \na^{\perp} \cdot u$ is described by the equation
\be \la{VOR}
\pa_t \omega +u\cdot \nabla \omega=-\nabla^{\perp}\rho \cdot \nabla \Phi.
\ee We take the $L^2$ inner product of \eqref{VOR} with $\omega$ and we estimate. 
The nonlinear term in $\omega$ vanishes due to the divergence-free condition obeyed by $u$. 
Elliptic estimates applied to the mean-free function $\rho$ yield
\be
\|\nabla \Phi\|_{L^\infty}\lesssim \|\rho\|_{L^4}\lesssim \|\nabla \rho\|_{L^2}.
\ee 
This gives the differential inequality 
\be
\begin{split}
\frac{1}{2}\frac{d}{dt}\|\omega\|_{L^2}^{2}
=
-\int \nabla^{\perp}\rho \cdot \nabla \Phi \omega\,dx
&\leq
\|\nabla \rho\|_{L^2}\|\nabla \Phi\|_{L^\infty}\|\omega\|_{L^2}
\lesssim  \|\na \rho\|_{L^2}^2 \|\omega\|_{L^2},
\end{split} 
\ee from which we conclude that
\be 
\fr{d}{dt} \|\omega\|_{L^2} \lesssim  \|\na \rho\|_{L^2}^2
\ee 
so that
\be\la{ol2}
\|\omega(t)\|_{L^2}\le \|\omega_0\|_{L^2}+C\int_0^t\|\na \rho(\tau)\|_{L^2}^2\,d\tau.
\ee
In view of the boundedness of $\rho$ in $L^2(0,T; H^1)$ obtained in Step 3, we conclude that $\omega$ belongs to $L^{\infty}(0,T; L^2)$. Since $\|\na u||_{L^2}\lesssim \|\omega\|_{L^2}$, Step 4 is complete.

\textbf{Step 5. A priori $H^2$ uniform-in-time bounds for the ionic concentrations.} Applying $\Delta$ to \eqref{1}, multiplying by $\Delta c_i$ and integrating by parts, we find that the $L^2$ norm of $\Delta c_i$ evolves according to 
\be \la{st6t1}
\frac{1}{2}\frac{d}{dt}\|\Delta c_i\|_{L^2}^2
+ D_i \|\na \Delta c_i\|_{L^2}^2
= - \int \Delta(u\cdot \nabla c_i) \Delta c_i \,dx
-D_i z_i\int  \Delta\nabla \cdot(c_i \nabla \Phi) \Delta c_i\,dx.
\ee For the nonlinear term in $u$, we integrate by parts and estimate using Sobolev and interpolation estimates. We obtain
\be 
\beg{aligned}
&\left|\int \Delta(u\cdot \nabla c_i) \Delta c_i \,dx \right|
= \left|\int \na ( u \cdot \na c_i) \cdot \na \Delta c_i \,dx \right|
\\&\le (\|\na u\|_{L^2}\|\na c_i\|_{L^{\infty}}+\|u\|_{L^4} \|\na\na c_i\|_{L^4}) \|\na \Delta c_i\|_{L^2}
\\&\lesssim\|\omega\|_{L^2}( \|\Delta c_i\|_{L^4}+\|\D c_i\|_{L^2}^{1/2}\|\na\D c_i\|_{L^2}^{1/2} )\|\na \Delta c_i\|_{L^2}\\
&\lesssim\|\omega\|_{L^2} \|\Delta c_i\|_{L^2}^{\fr{1}{2}} \|\na \Delta c_i\|_{L^2}^{\fr{3}{2}}
\\&\le \fr{D_i}{8} \|\na \Delta c_i\|_{L^2}^2 + C\|\omega\|_{L^2}^4 \|\Delta c_i\|_{L^2}^2.
\end{aligned}
\ee As for the second nonlinear term, we bound
\be \la{st61t1}
\beg{aligned}
&\left|-D_i z_i\int  \Delta\nabla \cdot(c_i \nabla \Phi) \Delta c_i\,dx \right|
= \left|D_i z_i \int \Delta (c_i \na \Phi) \cdot \na \Delta c_i \,dx\right|
\\&\lesssim\|\na \Delta c_i\|_{L^2} \left[\|\Delta c_i\|_{L^2} \|\na \Phi\|_{L^{\infty}} + \|c_i\|_{L^{4}}\|\na \rho\|_{L^4} \right]
\\&\le\fr{D_i}{8} \|\na \Delta c_i\|_{L^2}^2 + C\|\na \rho\|_{L^2}^2\|\Delta c_i\|_{L^2}^2 + C\|c_i\|_{L^4}^2 \|\Delta \rho\|_{L^2}^2
\end{aligned} 
\ee using fractional product estimates \cite{KP}. Putting \eqref{st6t1}--\eqref{st61t1} together, we obtain the differential inequality 
\be 
\fr{d}{dt} \|\Delta c_i\|_{L^2}^2 + D_i \|\na \Delta c_i\|_{L^2}^2
 \lesssim \left[\|\omega\|_{L^2}^4 + \|\na \rho\|_{L^2}^2\right] \|\Delta c_i\|_{L^2}^2
+ \|c_i\|_{L^4}^2 \|\Delta \rho\|_{L^2}^2.
\ee 
Using the triangle inequality, we bound 
\be 
\|\Delta \rho\|_{L^2}^2  \lesssim \sum\limits_{j=1}^{n} \|\Delta c_j\|_{L^2}^2.
\ee We sum over all indices $i \in \left\{1, ..., n \right\}$, and we obtain 
\be 
\beg{aligned}
&\fr{d}{dt} \left[\sum\limits_{i=1}^{n} \|\Delta c_i\|_{L^2}^2\right] + \sum\limits_{i=1}^{n} D_i \|\na \Delta c_i\|_{L^2}^2
\\&\lesssim  \left[\|\omega\|_{L^2}^4 + \|\na \rho\|_{L^2}^2\right] \left[\sum\limits_{i=1}^{n} \|\Delta c_i\|_{L^2}^2\right]
+ \left[\sum\limits_{i=1}^{n} \|c_i\|_{L^4}^2\right] \left[\sum\limits_{i=1}^{n} \|\Delta c_i\|_{L^2}^2 \right],
\end{aligned}
\ee 
and from Gronwall's inequality, we obtain
\be\la{ch2}
\begin{aligned}
&\sum_{i=1}^n\|\D c_i(t)\|_{L^2}^2+\sum_{i=1}^n\int_0^t D_i\|\na \D c_i(\tau)\|_{L^2}^2\,d\tau\\
\le &\sum_{i=1}^{n} \|\D c_i(0)\|_{L^2}^2\exp\left(C\int_0^T\|\omega(\tau)\|_{L^2}^4+\|\na\rho(\tau)\|_{L^2}^2+\sum_{i=1}^n\|c_i(\tau)\|_{L^4}^2\,d\tau\right).
\end{aligned}
\ee
From the above bounds, \eqref{Gamma} and \eqref{ol2}, we conclude that 
\be 
c_i \in L^{\infty}(0,T; H^2) \cap L^2(0,T;H^3)
\ee for all $i \in \left\{1, ..., n \right\}.$

\textbf{Step 6. A priori $L^\infty$ uniform-in-time bounds for the vorticity $\omega$.} We multiply the vorticity equation \eqref{VOR} by $\omega |\omega|^{p-2}$, integrate in the space variable over $\TT^2$, and estimate the forcing term in the resulting equation. We obtain 
\be
\fr{1}{p} \frac{d}{dt}\|\omega\|_{L^p}^{p}
\le \left| \int (\nabla^{\perp}\rho \cdot \nabla \Phi)\omega |\omega|^{p-2}\,dx\right|
\le \|\omega\|_{L^p}^{p-1}\|\nabla \rho\|_{L^p}\|\nabla \Phi\|_{L^\infty}
\ee by H\"older's inequality. In view of Morrey's inequality
\be 
\|\na \rho\|_{L^p} \le (4\pi^2)^{\fr{1}{p}} \|\na \rho\|_{L^{\infty}} \le  C_0(4\pi^2)^{\fr{1}{p}} \|\Delta \rho\|_{L^4} \le C_0(4\pi^2)^{\fr{1}{p}} \|\na \Delta \rho\|_{L^2}
\ee where $C_0$ is a constant independent of $p$. Therefore, we have
\be 
\fr{d}{dt} \|\omega\|_{L^p} \le  C_0(4\pi^2)^{\fr{1}{p}} \|\na \Phi\|_{L^{\infty}} \|\na \Delta \rho\|_{L^2}.
\ee We integrate in time from $0$ to $t$ and let $p \rightarrow \infty$. We obtain 
\be 
\|\omega(t)\|_{L^{\infty}} \le \|\omega_0\|_{L^{\infty}} + C_0 \int_{0}^{t} \|\na \Phi(s)\|_{L^{\infty}} \|\na \Delta \rho(s) \|_{L^2} ds
\ee for all $t \ge 0$. In view of \eqref{ch2}, we conclude that $\omega \in L^{\infty}(0,T; L^{\infty})$.

\textbf{Step 7. A priori $H^3$ uniform-in-time bounds for the velocity $u$.} 
For this step we outline the argument provided in \cite{IS}. We have 
\be \la{h3u}
\beg{aligned}
\fr{d}{dt} \|\na \Delta u\|_{L^2}^2
&\lesssim \|\omega\|_{L^{\infty}} \left(1 + \log (1 + \|\na \Delta u\|_{L^2}) \right) \|\na \Delta u\|_{L^2} 
\\&+ \|\rho\|_{L^{\infty}} \|\Delta \rho\|_{L^2}
+ \|\na \Phi\|_{L^{\infty}} \|\na \Delta \rho\|_{L^2},
\end{aligned}
\ee from which we conclude that
\be \la{h3u'}
\beg{aligned}
&\log\left(1 + \log(1 + \|\na \Delta u(t)\|_{L^2}) \right)
\le \log\left(1 + \log(1 + \|\na \Delta u_0\|_{L^2}) \right)
\\&+ C\int_{0}^{t} \left(\|\omega(\tau)\|_{L^{\infty}} + \|\rho(\tau)\|_{L^\infty}\|\Delta\rho(\tau)\|_{L^2}+\|\na\Phi(\tau)\|_{L^\infty}\|\na\D\rho(\tau)\|_{L^2} \right) d\tau
\end{aligned} 
\ee  for all $t \in [0,T]$.
This gives the double exponential bound
\be\la{uh3}
\|\na \D u(t)\|_{L^2}\lesssim (1 + \|\na \Delta u_0\|_{L^2}) \exp  \left(\exp \left(\int_{0}^{t} \left[\|\omega(\tau)\|_{L^{\infty}} + \|\na \Delta \rho (\tau)\|_{L^2}^2 \right] d\tau \right)\right)
\ee after applications of the Poincar\'e inequality. 
The regularity established in \eqref{ch2} allows us to conclude that $u \in L^{\infty}(0,T; H^3)$. We refer the reader to \cite[(3.23)--(3.25)]{IS} for the details involved in the preceding computations.

\textbf{Step 8. A priori $H^m$ uniform-in-time bounds for the ionic concentrations and velocity.} We take the scalar product in $L^2$ of each ionic concentration equation \eqref{1} with $\l^{2m} c_i$ and the vorticity equation \eqref{VOR} with $\l^{2m-2} \omega$. We add the resulting energy equations and obtain 
\be 
\beg{aligned} \la{st8t1}
&\frac{1}{2}\frac{d}{dt} \left[\|\Lambda^{m-1}\omega\|_{L^2}^{2} + \sum\limits_{i=1}^{n} \|\Lambda^m c_i\|_{L^2}^2\right]
+ \sum\limits_{i=1}^{n}  D_i \|\Lambda^{m+1}c_i\|_{L^2}^2
\\&= -\int \Lambda^{m-1}(u\cdot \nabla \omega) \Lambda^{m-1} \omega\,dx
- \int \Lambda^{m-1}(\nabla^{\perp}\rho \cdot \nabla \Phi) \Lambda^{m-1} \omega\,dx
\\&-\sum\limits_{i=1}^{n} \int \Lambda^{m} (u\cdot \nabla c_i) \Lambda^m c_i\,dx
+ \sum\limits_{i=1}^{n} D_iz_i \int \Lambda^{m} \nabla \cdot (c_i \nabla \Phi) \Lambda^m c_i\,dx.
\end{aligned}
\ee
We estimate 
\be 
\beg{aligned}
&\left|\int \Lambda^{m-1}(u\cdot \nabla \omega)\Lambda^{m-1} \omega\,dx \right|
= \left|\int [\Lambda^{m-1}, u\cdot \nabla ] \omega \l^{m-1} \omega \,dx \right|
\\&\le \|[\Lambda^{m-1}, u\cdot \nabla ]\omega\|_{L^2} \|\Lambda^{m-1} \omega\|_{L^2}
\lesssim (\|\nabla u\|_{L^\infty}\|\l^{m-1}\omega\|_{L^2}+\|\l^{m-1}u\|_{L^4}\|\na\omega\|_{L^4} )\|\l^{m-1} \omega\|_{L^2}\\
&\lesssim \|\na \Delta u\|_{L^2} \|\l^{m-1} \omega\|_{L^2}^2
\end{aligned}
\ee using standard commutator estimates (see \cite{KPV}). We bound
\be 
\beg{aligned}
&\left| \int \Lambda^{m-1}(\nabla^{\perp}\rho \cdot \nabla \Phi) \Lambda^{m-1} \omega\,dx \right|
\\&\lesssim \|\Lambda^{m-1} \omega\|_{L^2}
\left(\|\l^{m} \rho\|_{L^2}\|\nabla \Phi\|_{L^\infty} + \|\l^{m-1} \na \Phi\|_{L^4}\|\nabla \rho\|_{L^4}\right)
\\&\le \tilde a(t)\left(\|\l^{m-1}\omega\|_{L^2}^2+\sum_{i=1}^n\|\l^mc_i\|_{L^2}^2\right)
\end{aligned} 
\ee by applying fractional product estimates and continuous Sobolev embeddings. Above the coefficient $\tilde a(t)$ is given by
\be
\tilde a(t)=C(\|\na\Phi(t)\|_{L^\infty}+\|\na\rho(t)\|_{L^4}).
\ee

Using the fact that the fluid velocity $u$ is divergence-free, integrating by parts, and applying fractional product, interpolation, and Poincar\'e inequalities, we estimate 
\be 
\beg{aligned}
&\left| \sum\limits_{i=1}^{n} \int \Lambda^{m} (u\cdot \nabla c_i) \Lambda^m c_i\,dx \right|
=\left| \sum\limits_{i=1}^{n} \int \Lambda^{m} (u c_i)\cdot \na \Lambda^m c_i\,dx \right|
\\&\lesssim \sum\limits_{i=1}^{n} \|\l^{m+1} c_i\|_{L^2} \left(\|\l^{m} u\|_{L^2} \|c_i\|_{L^{\infty}} + \|u\|_{L^{\infty}} \|\l^{m} c_i\|_{L^2} \right)
\\&\lesssim \sum\limits_{i=1}^{n} \|\l^{m+1} c_i\|_{L^2} \left(\|\l^{m-1} \omega\|_{L^2} \|c_i\|_{L^{\infty}} + \|\Delta u\|_{L^{2}} \|\l^{m} c_i\|_{L^2} \right)
\\&\le \sum\limits_{i=1}^{n} \fr{D_i}{8} \|\l^{m+1}c_i\|_{L^2}^2 + C\left(\sum\limits_{i=1}^{n} \|c_i\|_{L^{\infty}}^2 \right) \|\l^{m-1} \omega\|_{L^2}^2 + C\|\Delta u\|_{L^2}^2 \left(\sum\limits_{i=1}^{n} \|\l^{m} c_i\|_{L^2}^2 \right)
\\&\le \sum\limits_{i=1}^{n} \fr{D_i}{8} \|\l^{m+1}c_i\|_{L^2}^2 + C\left(\|\Delta u\|_{L^2}^2 + \sum\limits_{i=1}^{n} \|c_i\|_{L^{\infty}}^2 \right) \left(\|\l^{m-1} \omega\|_{L^2}^2 + \sum\limits_{i=1}^{n} \|\l^{m} c_i\|_{L^2}^2 \right).
\end{aligned}
\ee 
Using elliptic regularity, we have
\be \la{st81t1}
\beg{aligned}
&\left|\sum\limits_{i=1}^{n} D_iz_i \int \Lambda^{m} \nabla \cdot (c_i \nabla \Phi) \Lambda^m c_i\,dx \right|
= \left|\sum\limits_{i=1}^{n} D_iz_i \int \Lambda^{m}  (c_i \nabla \Phi) \cdot \na \Lambda^m c_i\,dx \right|
\\&\lesssim \sum\limits_{i=1}^{n} \|\l^{m+1} c_i\|_{L^2} \left(\|\l^{m} c_i\|_{L^2} \|\na \Phi\|_{L^{\infty}} + \|c_i\|_{L^4} \|\l^{m} \na \Phi\|_{L^4} \right)
\\&\lesssim\sum\limits_{i=1}^{n} \|\l^{m+1} c_i\|_{L^2} \left(\|\l^{m} c_i\|_{L^2} \|\na \rho\|_{L^{2}} + \|c_i\|_{L^4} \|\l^{m} \rho\|_{L^2} \right)
\\&\le \sum\limits_{i=1}^{n} \fr{D_i}{8} \|\l^{m+1} c_i\|_{L^2}^2 + C\left(\|\na \rho\|_{L^2}^2 + \sum\limits_{i=1}^{n} \|c_i\|_{L^4}^2 \right)\left(\sum\limits_{i=1}^{n} \|\l^{m} c_i\|_{L^2}^2 \right).
\end{aligned}
\ee Putting \eqref{st8t1}--\eqref{st81t1} together, we obtain the differential inequality
\be \la{mb'}
\frac{d}{dt} \left[\|\Lambda^{m-1}\omega\|_{L^2}^{2} + \sum\limits_{i=1}^{n} \|\Lambda^m c_i\|_{L^2}^2\right]
\lesssim a(t) \left[\|\Lambda^{m-1}\omega\|_{L^2}^{2} + \sum\limits_{i=1}^{n} \|\Lambda^m c_i\|_{L^2}^2\right]
\ee where 
\be 
a(t) = \tilde a(t)+\|\na\D u(t)\|_{L^2}+\| \Delta u(t)\|_{L^2}^2 + \|\nabla \rho(t)\|_{L^2}^2  + \sum_{i=1}^{n} \|c_i(t)\|_{L^{\infty}}^2
\ee 
is integrable over the time interval $[0,T]$ due to \eqref{ch2} and \eqref{uh3}. Thus, applying Gronwall's inequality to \eqref{mb'}, we obtain
\be
\|\Lambda^{m-1}\omega(t)\|_{L^2}^{2} + \sum\limits_{i=1}^{n} \|\Lambda^m c_i(t)\|_{L^2}^2\le \left[\|\Lambda^{m-1}\omega_0\|_{L^2}^{2} + \sum\limits_{i=1}^{n} \|\Lambda^m c_i(0)\|_{L^2}^2\right]\exp\left(C\int_0^t a(\tau)\,d\tau\right)
\ee
This finishes the proof of Step 8.

\textbf{Step 9. Extension of the local solution.} The local solution can be extended to the time interval $[0,T]$, a fact that follows from the uniform-in-time boundedness in $H^m$ obtained in Step 8. This ends the proof of Theorem~\ref{T1}.

\section{Gevrey Regularity of the Nernst-Planck-Euler system}\la{S3}

In this section, we address the propagation of the Gevrey regularity for the NPE system:

\beg{thm} \la{T2} (Global analyticity of NPE) Let $T>0$ be arbitrary and $m > 4$. Assume that the initial ionic concentrations $c_i(0)$ are nonnegative and real-analytic, and the initial velocity $u_0$ is mean-zero, divergence-free, and real-analytic.  Then the NPE system  described by \eqref{1}, \eqref{4} and \eqref{2} has a unique solution $(u, c_1, ..., c_n)$ on $[0,T]$ such that for each time $t \in (0,T)$, the functions $u, c_1,..., c_n$ are real analytic in the spatial variable with uniform radius of analyticity $\tau(t)$, depending on the $H^m$ norm of the solution $(u, c_1, ..., c_n)$ up to time $t$.
\end{thm}

In order to prove Theorem \ref{T2}, we start by taking the scalar products in $\mathcal{D}(e^{\tau \Lambda})$ of the equation obeyed by the ionic concentration $c_i$ with $\l^{2m}c_i$ and of the equation obeyed by $\omega = \na^{\perp} \cdot u$ with $\l^{2m-2}\omega$. We add the resulting energy equalities, and we obtain
\be \la{energy}
\beg{aligned} 
&\fr{1}{2} \fr{d}{dt} \left[\sum\limits_{i=1}^{n} \|e^{\tau \l} \l^m c_i\|_{L^2}^2 + \|e^{\tau \l} \l^{m-1} \omega\|_{L^2}^2 \right]
+ \sum\limits_{i=1}^{n} D_i \|e^{\tau \l} \l^{m+1} c_i\|_{L^2}^2 
\\&= \tau'(t) \sum\limits_{i=1}^{n} \|e^{\tau \l} \l^{m + \fr{1}{2}} c_i\|_{L^2}^2 
+ \tau'(t) \|e^{\tau \l} \l^{m - \fr{1}{2}} \omega\|_{L^2}^2 
- \sum\limits_{i=1}^{n} (e^{\tau \l} \l^m (u \cdot \na (c_i - \bar{c_i}) ), e^{\tau \l} \l^m(c_i - \bar{c}_i))_{L^2}
\\&\quad+ \sum\limits_{i=1}^{n} D_i z_i (e^{\tau \l} \l^{m} \na \cdot [(c_i - \bar{c}_i) \na \Phi], e^{\tau \l} \l^m (c_i - \bar{c}_i))_{L^2}
+ \sum\limits_{i=1}^{n} D_i z_i \bar{c}_i (e^{\tau \l}\l^m \Delta \Phi, e^{\tau \l} \l^m (c_i - \bar{c}_i))_{L^2} 
\\&\quad\quad- (u \cdot \na \omega, \l^{2m-2}e^{2\tau \l} \omega)_{L^2} 
- (\na^{\perp} \cdot [\rho \na \Phi], \l^{2m-2} e^{2\tau \l} \omega)_{L^2}.
\end{aligned}
\ee

In order to control the nonlinear terms, we need the following lemmas.

\beg{lem} \cite{KV} The following estimate
\be \la{nonlin1}
\beg{aligned}  
&|(u \cdot \na \omega, \l^{2m-2} e^{2\tau \l} \omega)_{L^2}| 
\\&\quad\quad\lesssim  \left(\tau \|\na u\|_{L^{\infty}} + \tau^2 \|\l^{m-1}\omega\|_{L^2} + \tau^2 \|e^{\tau \l} \l^{m-1} \omega\|_{L^2} \right) \|e^{\tau \l} \l^{m - \fr{1}{2}} \omega\|_{L^2}^2 
\\&\quad\quad\quad\quad \left(\|\na u\|_{L^{\infty}} \|\l^{m-1}e^{\tau \l}\omega\|_{L^2} + (1+\tau)\|\l^{m-1} \omega\|_{L^2}^2 \right) \|\l^{m-1} e^{\tau \l} \omega\|_{L^2}
\end{aligned}  
\ee holds for any $m > 4$.
\end{lem}

\beg{lem} \la{L2} Let $i \in \left\{1, ..., n \right\}$. The following estimate
\be \la{nonlin2}
\beg{aligned} 
&|(e^{\tau \l} \l^{m} (u \cdot \na (c_i - \bar{c}_i )), e^{\tau \l} \l^m (c_i - \bar{c}_i ))_{L^2}| 
\\&\quad\quad\leq  \fr{D_i}{8} \|\l^{m+1} e^{\tau \l} c_i\|_{L^2}^2
+ C\|\l^{1+\epsilon} u\|_{L^2}^2 \|\l^m e^{\tau \l} c_i\|_{L^2}^2
+ C \|\l^{m-1} u\|_{L^2}^2 \|\l^{2+\epsilon} e^{\tau \l} c_i\|_{L^2}^2 
\\&\quad\quad\quad+ C \tau^2 \|\l^{2+\epsilon} e^{\tau \l} u\|_{L^2}^2 \|\l^{m} e^{\tau \l} c_i\|_{L^2}^2
+ C\tau^2 \|\l^{m} e^{\tau \l} u\|_{L^2}^2 \|\l^{2+\epsilon} e^{\tau \l} c_i\|_{L^2}^2.
\end{aligned} \ee holds for any $m > 0$ and $\epsilon>0$.
Consequently, we have 
\be \la{nonlin22}
\beg{aligned} 
&|(e^{\tau \l} \l^{m} (u \cdot \na (c_i - \bar{c}_i )), e^{\tau \l} \l^m (c_i - \bar{c}_i ))_{L^2}| 
\\&\quad\quad\leq  \fr{D_i}{8} \|\l^{m+1} e^{\tau \l} c_i\|_{L^2}^2
+ C \|\l^{m-1} u\|_{L^2}^2 \|\l^{m} e^{\tau \l} c_i\|_{L^2}^2 
+ C \tau^2  \|\l^{m} e^{\tau \l} c_i\|_{L^2}^2 \|\l^{m-\fr{1}{2}} e^{\tau \l} \omega\|_{L^2}^2 
\end{aligned}
\ee
for any $m > 2$.
\end{lem}

\textbf{Proof.} We set 
\be \la{ufourier}
u = \sum\limits_{j \in \ZZ^2 \setminus \left\{0\right\}} u_j e^{ij \cdot x},
\ee
\be \la{ustarfourier}
u^* = \sum\limits_{j \in \ZZ^2 \setminus \left\{0\right\}} u_j^* e^{ij \cdot x}, u_j^* = e^{\tau|j|} u_j,
\ee
\be \la{cfourier}
c_i - \bar{c}_i=   \sum\limits_{j \in \ZZ^2 \setminus \left\{0\right\}} (c_i)_j e^{ij \cdot x},
\ee and
\be \la{cstarfourier}
(c_i - \bar{c}_i)^* = \sum\limits_{j \in \ZZ^2 \setminus \left\{0\right\}} (c_i)_j^* e^{ij \cdot x}, (c_i)_j^* = e^{\tau|j|} (c_i)_j.
\ee 
Then 
\be 
\left|(e^{\tau \l} \l^{m} (u \cdot \na (c_i - \bar{c}_i )), e^{\tau \l} \l^m (c_i - \bar{c}_i ))_{L^2}\right| 
= 4\pi^2 \left|\sum\limits_{j+k+l = 0} (u_j \cdot k) (c_i)_{k} (c_i)_{l} e^{2\tau |l|} |l|^{2m} \right|.
\ee Since $j+k + l = 0$, we have $|l| \le |k| + |j|$, which implies that $e^{\tau |l|} \le e^{\tau |k|}e^{\tau |j|}$, and so 
\be 
\left|(e^{\tau \l} \l^{m} (u \cdot \na (c_i - \bar{c}_i )), e^{\tau \l} \l^m (c_i - \bar{c}_i ))_{L^2}\right| 
\le 4\pi^2 \sum\limits_{j+k+l =0} e^{\tau |j|} |u_j \cdot k| |(c_i)_{k}^*| |(c_i)_{l}^*|  |l|^{2m} .
\ee Using the estimate $e^x \le e + xe^x$, that holds for any $x \ge 0$, and the estimate $|l|^{m-1}  \lesssim (|k|^{m-1} + |j|^{m-1})$, that holds for all $j, k, l \in \ZZ^2$ with $j+k+l = 0$, we bound
\be 
\beg{aligned}
&\left|(e^{\tau \l} \l^{m} (u \cdot \na (c_i - \bar{c}_i )), e^{\tau \l} \l^m (c_i - \bar{c}_i ))_{L^2}\right|  
\\&\quad\quad\lesssim \sum\limits_{j+k+l = 0} (e + \tau |j| e^{\tau |j|}) |u_j| |(c_i)_{k}^*| |(c_i)_{l}^*| |k| |l|^{m+1} (|k|^{m-1} + |j|^{m-1}).
\end{aligned}
\ee Using Young's convolution inequality and Plancherel identity, we obtain 
\be
\beg{aligned}
&\left|(e^{\tau \l} \l^{m} (u \cdot \na (c_i - \bar{c}_i )), e^{\tau \l} \l^m (c_i - \bar{c}_i ))_{L^2}\right| 
\\&\quad\quad\lesssim \|u_j\|_{\ell^1} \||k|^m (c_i)_{k}^*\|_{\ell^2} \||l|^{m+1} (c_i)_l^*\|_{\ell^2} 
+ \||j|^{m-1}u_j\|_{\ell^2} \||k| (c_i)_k^*\|_{\ell^1} \||l|^{m+1}(c_i)_l^*\|_{\ell^2} 
\\&\quad\quad\quad+ \tau \||j|u_j^*\|_{\ell^1}\||k|^m (c_i)_{k}^*\|_{\ell^2} \||l|^{m+1} (c_i)_l^*\|_{\ell^2}
+ \tau \||j|^m u_j^*\|_{\ell^2}  \||k| (c_i)_k^*\|_{\ell^1} \||l|^{m+1}(c_i)_l^*\|_{\ell^2} 
\\&\quad\quad\lesssim \|\l^{1+\epsilon} u\|_{L^2} \|\l^{m} c_i^*\|_{L^2} \|\l^{m+1} c_i^*\|_{L^2}
+ \|\l^{m-1}u\|_{L^2} \|\l^{2+\epsilon} c_i^*\|_{L^2} \|\l^{m+1} c_i^*\|_{L^2} 
\\&\quad\quad\quad+ \tau \|\l^{2+\epsilon} u^*\|_{L^2} \|\l^{m} c_i^*\|_{L^2} \|\l^{m+1} c_i^*\|_{L^2}
+ \tau \|\l^{m} u^*\|_{L^2} \|\l^{2+\epsilon} c_i^*\|_{L^2} \|\l^{m+1} c_i^*\|_{L^2}
\end{aligned}
\ee which gives the desired estimate \eqref{nonlin2}. The bound \eqref{nonlin22} is a consequence of Young's inequality and the bound
\be 
\|\l^{m} e^{\tau \l} u\|_{L^2} = \|\l \l^{m-1} e^{\tau \l} u\|_{L^2} \lesssim \|\na \l^{m-1} e^{\tau \l} u\|_{L^2} 
\lesssim \|\l^{m-1} e^{\tau \l} \na^{\perp} \cdot u \|_{L^2} \lesssim \|\l^{m-\fr{1}{2}} e^{\tau \l} \omega\|_{L^2}
\ee that holds in view of the divergence-free condition obeyed by the velocity $u$. This ends the proof of Lemma~\ref{L2}.

\beg{lem} \la{L3}  Let $i \in \left\{1, ..., n \right\}$. The following estimate 
\be \la{nonlin3}
\beg{aligned} 
&D_i |z_i||(e^{\tau \l} \l^{m} \na \cdot [(c_i - \bar{c}_i) \na \Phi], e^{\tau \l} \l^m (c_i - \bar{c}_i))_{L^2}| 
\\&\quad\quad\le \fr{D_i}{8} \|\l^{m+1} e^{\tau \l} c_i\|_{L^2}^2 
+ C \|\l^{m-1} \rho\|_{L^2}^2 \|\l^{1+\epsilon} e^{\tau \l} c_i\|_{L^2}^2
+ C  \|\l^{\epsilon} \rho \|_{L^2}^2 \|\l^{m} e^{\tau \l} c_i\|_{L^2}^2 
\\&\quad\quad\quad+ C\tau^2 \|\l^{1+\epsilon} e^{\tau \l} c_i\|_{L^2}^2 \|\l^{m} e^{\tau \l} \rho\|_{L^2}^2
+ C \tau^2 \|\l^{1 + \epsilon} e^{\tau \l} \rho \|_{L^2}^2 \|\l^{m} e^{\tau \l} c_i\|_{L^2}^2
\end{aligned} \ee holds for any $m > 0$ and $\epsilon>0$. Consequently, we have
\be  \la{nonlin32}
\beg{aligned} 
&D_i |z_i| |(e^{\tau \l} \l^{m} \na \cdot [(c_i - \bar{c}_i) \na \Phi], e^{\tau \l} \l^m (c_i - \bar{c}_i))_{L^2}| 
\\&\quad\quad\le \fr{D_i}{8} \|\l^{m+1} e^{\tau \l} c_i\|_{L^2}^2 
+ C \|\l^{m-1} \rho\|_{L^2}^2 \|\l^{m} e^{\tau \l} c_i\|_{L^2}^2 
+ C\tau^2  \|\l^{m} e^{\tau \l} \rho\|_{L^2}^2 \|\l^{m+ \fr{1}{2}} e^{\tau \l} c_i\|_{L^2}^2
\end{aligned} \ee for any $m > 1$.
\end{lem} 

\textbf{Proof.} We set $c_i$ and $c_i^*$ as in \eqref{cfourier} and \eqref{cstarfourier} respectively. We write the Fourier series of $\rho$ as 
\be 
\rho = \sum\limits_{k \in \ZZ^2 \setminus \left\{0\right\}} \rho_k e^{ik \cdot x},
\ee and we set 
\be 
\rho^* = \sum\limits_{k \in \ZZ^2 \setminus \left\{0\right\}} \rho_k^* e^{ik \cdot x}, \rho_k^* = e^{\tau|k|} \rho_k.
\ee
The gradient of the potential $\Phi$ has a Fourier expansion given by 
\be 
\na \Phi = \na (-\Delta)^{-1} \rho = i \sum\limits_{k \in \ZZ^2 \setminus \left\{0\right\}} \fr{k}{|k|^2} \rho_k e^{ik \cdot x}.
\ee We have 
\be \la{lem30}
(e^{\tau \l} \l^{m} \na \cdot [(c_i - \bar{c}_i) \na \Phi], e^{\tau \l} \l^m (c_i - \bar{c}_i))_{L^2}
= -4\pi^2 \sum\limits_{j + k + l = 0} (c_i)_{j} \rho_k (c_i)_l \fr{l \cdot k}{|k|^2}|l|^{2m} e^{2\tau |l|},
\ee hence
\be \la{lem31}
\beg{aligned} 
&|(e^{\tau \l} \l^{m} \na \cdot [(c_i - \bar{c}_i) \na \Phi], e^{\tau \l} \l^m (c_i - \bar{c}_i))_{L^2}| 
\\&\quad\quad \lesssim \sum\limits_{j+k + l = 0} |(c_i)_j^*| |(c_i)_l^*| (e + \tau |k| e^{\tau |k|}) |\rho_k||k|^{-1} |l|^{m+1} (|k|^m + |j|^m)
\end{aligned} \ee using the triangle inequality $|l| \le |k| + |j|$. In view of Young's convolution inequality, we estimate the following four terms
\be 
\beg{aligned}
&D_i |z_i| \sum\limits_{j+k+l = 0} e |(c_i)_j^*| |(c_i)_l^*| |\rho_k| |l|^{m+1} |k|^{m-1} 
\lesssim \||l|^{m+1} (c_i)_l^*\|_{\ell^2} \||k|^{m-1} \rho_k \|_{\ell^2} \|(c_i)_j^*\|_{\ell^1}
\\&\quad\quad\le \fr{D_i}{32} \|\l^{m+1} e^{\tau \l} c_i\|_{L^2}^2 + C \|\l^{1+\epsilon} e^{\tau \l} c_i\|_{L^2}^2 \|\l^{m-1} \rho\|_{L^2}^2,
\end{aligned} \ee
\be 
\beg{aligned}
&D_i |z_i| \sum\limits_{j+k+l = 0} e |(c_i)_j^*| |(c_i)_l^*||\rho_k| |l|^{m+1} |k|^{-1}|j|^m 
\lesssim \||l|^{m+1} (c_i)_l^*\|_{\ell^2} \||j|^{m} (c_i)_j^* \|_{\ell^2} \||k|^{-1} \rho_k\|_{\ell^1}
\\&\quad\quad\le \fr{D_i}{32} \|\l^{m+1} e^{\tau \l} c_i\|_{L^2}^2 + C\|\l^{\epsilon} \rho \|_{L^2}^2 \|\l^{m} e^{\tau \l} c_i\|_{L^2}^2,
\end{aligned} \ee
\be 
\beg{aligned}
&D_i |z_i| \sum\limits_{j+k+l = 0} \tau |(c_i)_j^*| |(c_i)_l^*| |\rho_k^*| |l|^{m+1} |k|^{m} 
 \lesssim \tau \||l|^{m+1} (c_i)_l^*\|_{\ell^2} \||k|^{m} \rho_k^* \|_{\ell^2} \|(c_i)_j^*\|_{\ell^1}
\\&\quad\quad\le \fr{D_i}{32} \|\l^{m+1} e^{\tau \l} c_i\|_{L^2}^2 + C \tau^2\|\l^{1+\epsilon} e^{\tau \l} c_i\|_{L^2}^2 \|\l^{m} e^{\tau \l} \rho\|_{L^2}^2,
\end{aligned} \ee  and
\be \la{lem32}
\beg{aligned} 
&D_i |z_i| \sum\limits_{j+k+l = 0} \tau |(c_i)_j^*| |(c_i)_l^*| |\rho_k^*| |l|^{m+1} |j|^m 
 \lesssim \tau \||l|^{m+1} (c_i)_l^*\|_{\ell^2} \||j|^{m} (c_i)_j^* \|_{\ell^2} \| \rho_k^*\|_{\ell^1}
\\&\quad\quad\le \fr{D_i}{32} \|\l^{m+1} e^{\tau \l} c_i\|_{L^2}^2 + C\tau^2 \|\l^{1 + \epsilon} e^{\tau \l} \rho \|_{L^2}^2 \|\l^{m} e^{\tau \l} c_i\|_{L^2}^2.
\end{aligned} \ee
Putting \eqref{lem31}--\eqref{lem32} together, we obtain the desired estimate \eqref{nonlin3}. The bound \eqref{nonlin32} follows from \eqref{nonlin3} and continuous embeddings of Sobolev spaces. This completes the proof of Lemma \ref{L3}.

\beg{lem} \la{L4} Let $i \in \left\{1, ..., n \right\}$. The following estimate
\be \la{nonlin4}
\beg{aligned} 
&D_i |z_i| |\bar{c}_i(0)| |(e^{\tau \l}\l^m \Delta \Phi, e^{\tau \l} \l^m (c_i - \bar{c}_i))_{L^2}| 
\\&\quad\quad\le \fr{D_i}{8}\|\l^{m+1}e^{\tau\l}c_i\|_{L^2}^2+C|\bar{c}_i(0)|^2\|\l^{m}e^{\tau\l}\rho\|_{L^2}^2
\end{aligned} \ee holds for any $m > 0 $.
\end{lem}

\textbf{Proof.} Under the same settings of the previous lemmas, we have 
\be 
    \beg{aligned}
    &D_i |z_i| |\bar{c}_i(0)| |(e^{\tau \l}\l^m \Delta \Phi, e^{\tau \l} \l^m (c_i - \bar{c}_i))_{L^2}|
    \lesssim |\bar{c}_i(0)| \sum_{k\in\mathbb{Z}^2\setminus\{0\}}e^{2\tau|k|}|k|^{2m}|\rho_k||(c_i)_k|
    \\&\quad\quad\lesssim |\bar{c}_i(0)|  \||k|^{m+1}(c_i)_k^*\|_{\ell^2}\||k|^{m-1}\rho_k^*\|_{\ell^2}
    \le \fr{D_i}{8}\|\l^{m+1}e^{\tau\l}c_i\|_{L^2}^2+C |\bar{c}_i(0)|^2\|\l^{m-1}e^{\tau\l}\rho\|_{L^2}^2
    \\&\quad\quad\le \fr{D_i}{8}\|\l^{m+1}e^{\tau\l}c_i\|_{L^2}^2+C|\bar{c}_i(0)|^2\|\l^{m}e^{\tau\l}\rho\|_{L^2}^2.
    \end{aligned}
\ee 
Here we used the fact $\rho_{-k}$ coincides with the complex conjugate of $\rho_{k}$, which follows from the real-valuedness of $\rho$. This finishes the proof of Lemma \ref{L4}.

\beg{lem} \la{L5} The following estimate
\be \la{nonlin5}
\beg{aligned} 
&|(\na^{\perp} \cdot [\rho \na \Phi], \l^{2m-2} e^{2\tau \l} \omega)_{L^2}| 
\\&\quad\quad\lesssim \|\l^{m-1} \rho\|_{L^2} \|\l^{1+\epsilon} e^{\tau \l}\rho\|_{L^2}^2  + \|\l^{\epsilon} \rho\|_{L^2} \|\l^{m} e^{\tau \l}\rho \|_{L^2}^2  + \| \l^{m-1} e^{\tau \l} \omega\|_{L^2}^2 \|\l^{m-1} \rho\|_{L^2}
\\&\quad\quad\quad\quad+ \| \l^{m-1} e^{\tau \l} \omega\|_{L^2}^2  \|\l^{\epsilon} \rho\|_{L^2} + \|\l^{m} e^{\tau \l} \rho\|_{L^2}^2+ \tau^2 \|\l^{1+\epsilon} e^{\tau \l} \rho\|_{L^2}^2 \| \l^{m-1} e^{\tau \l} \omega\|_{L^2}^2 
\end{aligned} \ee holds for any $m \ge 1$ and any $\epsilon > 0$. Consequently, we have 
\be \la{nonlin51}
\beg{aligned} 
&|(\na^{\perp} \cdot [\rho \na \Phi], \l^{2m-2} e^{2\tau \l} \omega)_{L^2}| 
\\&\quad\quad\lesssim \|\l^{m-1} \rho\|_{L^2} \|\l^{m} e^{\tau \l}\rho\|_{L^2}^2   + \|\l^{m-1} \rho\|_{L^2} \| \l^{m-1} e^{\tau \l} \omega\|_{L^2}^2  + \|\l^{m} e^{\tau \l} \rho\|_{L^2}^2 
\\&\quad\quad\quad\quad+ \tau^2 \|\l^{m} e^{\tau \l} \rho\|_{L^2}^2 \| \l^{m-\fr{1}{2}} e^{\tau \l} \omega\|_{L^2}^2 
\end{aligned}
\ee for any $m > 1$. 
\end{lem}

\textbf{Proof.} We set $\rho$ and $\rho^*$ as in the previous lemmas, and we set 
\be 
\omega = \sum\limits_{l \in \ZZ^2 \setminus \left\{0\right\}} \omega_l e^{il \cdot x}, \omega^* = \sum\limits_{l \in \ZZ^2 \setminus \left\{0\right\}} \omega_l^* e^{il \cdot x}, \omega_l^* = e^{\tau |l|} \omega_l. 
\ee We have 
\be \la{lem50}
\beg{aligned} 
&\left|(e^{\tau \l} \l^{m-1} \na^{\perp} \cdot [\rho \na \Phi], e^{\tau \l} \l^{m-1} \omega)_{L^2} \right|
= 4\pi^2 \left| \sum\limits_{j + k + l = 0} \rho_{j} \rho_k \omega_l \fr{l^{\perp} \cdot k}{|k|^2}|l|^{2m-2} e^{2\tau |l|} \right|
\\&\quad\quad\le \sum\limits_{j + k + l = 0} |\rho_j| |\rho_k| |\omega_l| |k|^{-1} |l|^{2m-1} e^{2\tau |l|}
\\&\quad\quad\lesssim \sum\limits_{j + k + l = 0} |\rho_j^*| (e + \tau|k|e^{\tau |k|})|\rho_k| |\omega_l^*| |k|^{-1} |l|^{m-1} \left(|k|^m + |j|^m\right). 
\end{aligned}
\ee We estimate the following two sums 
\be 
\beg{aligned}
&\sum\limits_{j + k + l = 0} e |\rho_j^*|  |\rho_k| |\omega_l^*| |k|^{-1} |l|^{m-1} \left(|k|^m + |j|^m\right)
\\&\quad\quad\lesssim \| |l|^{m-1} \omega_l^*\|_{\ell^2} \||k|^{m-1} \rho_k\|_{\ell^2} \|\rho_j^*\|_{\ell^1}
+ \| |l|^{m-1} \omega_l^*\|_{\ell^2} \||j|^{m} \rho_j^*\|_{\ell^2} \||k|^{-1}\rho_k\|_{\ell^1}
\\&\quad\quad\lesssim \| \l^{m-1} e^{\tau \l} \omega\|_{L^2} \|\l^{m-1} \rho\|_{L^2} \|\l^{1+\epsilon} e^{\tau \l}\rho\|_{L^2}
+ \| \l^{m-1} e^{\tau \l} \omega\|_{L^2} \|\l^{m} e^{\tau \l}\rho \|_{L^2} \|\l^{\epsilon} \rho\|_{L^2}
\end{aligned}
\ee and 
\be \la{lem51}
\beg{aligned}
&\sum\limits_{j + k + l = 0} \tau |\rho_j^*| |k|e^{\tau |k|}|\rho_k| |\omega_l^*| |k|^{-1} |l|^{m-1} \left(|k|^m + |j|^m\right)
\\&\quad\quad\lesssim \tau \| |l|^{m-1} \omega_l^*\|_{\ell^2} \||k|^{m} \rho_k^*\|_{\ell^2} \|\rho_j^*\|_{\ell^1}
+ \tau \| |l|^{m-1} \omega_l^*\|_{\ell^2} \||j|^{m} \rho_j^*\|_{\ell^2} \|\rho_k^*\|_{\ell^1}
\\&\quad\quad\lesssim \tau \| \l^{m-1} e^{\tau \l} \omega\|_{L^2} \|\l^m e^{\tau \l} \rho\|_{L^2} \|\l^{1+\epsilon} e^{\tau \l} \rho\|_{L^2}
+ \tau \| \l^{m-1} e^{\tau \l} \omega\|_{L^2} \|\l^m e^{\tau \l} \rho\|_{L^2}  \|\l^{1+\epsilon} e^{\tau \l} \rho\|_{L^2}.
\end{aligned}
\ee Putting \eqref{lem50}--\eqref{lem51} together and applying Young's inequality, we obtain the bound \eqref{nonlin5} for any $m \ge 1$ and $\epsilon > 0$. The bound \eqref{nonlin51} is a direct consequence of \eqref{nonlin5}. This ends the proof of Lemma \ref{L5}. 

We put \eqref{energy}, \eqref{nonlin1}, \eqref{nonlin22}, \eqref{nonlin32}, \eqref{nonlin4}, and \eqref{nonlin51} together. Setting 
\be 
y(t) = \sum\limits_{i=1}^{n} \|e^{\tau \l} \l^m c_i\|_{L^2}^2 + \|e^{\tau \l} \l^{m-1} \omega\|_{L^2}^2,
\ee
\be 
\tilde{B}(t) = \tau'(t) + C\tau \|\na u\|_{L^{\infty}} + C\tau^2 \|\l^{m-1}\omega\|_{L^2} + C\tau^2 \|e^{\tau \l} \l^{m-1} \omega\|_{L^2}+ C\tau^2  \sum\limits_{i=1}^{n} \|\l^{m} e^{\tau \l} c_i\|_{L^2}^2,
\ee and
\be 
B(t) = \|\l^{m-1} \rho\|_{L^2} + \|\l^{m-1} \rho\|_{L^2}^2 + \chi_{\sum\limits_{i=1}^{n} c_i(0) \ne 0}  \left(\|\l^{m-1} u\|_{L^2}^2 + 1 \right) + \|\na u\|_{L^{\infty}} + \sum\limits_{i=1}^{n} |\bar{c}_i(0)|^2,
\ee where $ \chi_{\sum\limits_{i=1}^{n} c_i(0) \ne 0}$ is the characteristic function of the singleton set $\left\{\sum\limits_{i=1}^{n} c_i(0) \ne 0 \right\}$, we obtain the energy differential inequality 
\be 
\beg{aligned}
\fr{1}{2} \fr{d}{dt} y(t) &\le \tilde{B}(t)\left[\sum\limits_{i=1}^{n} \|e^{\tau \l} \l^{m + \fr{1}{2}} c_i\|_{L^2}^2 + \|e^{\tau \l} \l^{m - \fr{1}{2}} \omega\|_{L^2}^2  \right] 
\\&\quad\quad+ CB(t) y(t) + C(1+\tau) \|\l^{m-1} \omega\|_{L^2}^2 \sqrt{y(t)}
\end{aligned} 
\ee for any $m > 4$. Above the characteristic function is introduced to explicitly indicate terms that would be absent if we assumed the absence of ions so that the NPE system reduces to the Euler equations (c.f. Remark \ref{RM}).

If $\tau$ is such that $\tilde{B} (t) \le 0$ (so that in particular $\tau'\le 0$), then the latter inequality reduces to 
\be 
\fr{d}{dt} \sqrt{y(t)} \lesssim B(t)\sqrt{y(t)} + (1 + \tau(0)) \|\l^{m-1} \omega\|_{L^2}^2.
\ee Therefore, we obtain 
\be \la{sc}
\sqrt{y(t)} \le g(t) \left\{ \sqrt{y(0)} + C(1 + \tau(0)) \int_{0}^{t} \|\l^{m-1} \omega(s)\|_{L^2}^2 g(s)^{-1} ds \right\} :=A(t)
\ee for any $t \in [0,T]$, where 
\be 
g(t) = C\int_{0}^{t} B(s) ds.
\ee In light of \eqref{sc}, a sufficient condition for $\tilde{B}(t)$ to be nonpositive is that
\be \label{e.w11021}
\tau'(t) + C\tau \|\na u\|_{L^{\infty}} + C\tau^2 \|\l^{m-1}\omega\|_{L^2} + C\tau^2 \tilde{A}(t) \le 0
\ee  where
\be 
\tilde{A} (t) = A(t) +  \chi_{\sum\limits_{i=1}^{n} c_i(0) \ne 0}  A(t)^2.
\ee Hence, it suffices to choose
	\begin{equation} \la{radius}
	\tau(t)=\frac{1}{g(t) \left(\frac{1}{\tau(0)}
		+C\int_{0}^{t} (\|\Lambda^{m-1}\omega(s)\|_{L^2}+\tilde{A}(s))g^{-1}(s)\,ds
		\right)}.
	\end{equation}
This completes the proof of Theorem \ref{T2}.

\beg{rem}\la{RM} If the initial ionic concentrations are taken to be zero, then their spatial averages vanish at all positive times, and so they are identically zero on $\TT^2 \times [0, \infty)$. In this case, the NPE problem reduces to the two-dimensional periodic non-forced Euler system and yields the spatial analyticity of its solution. In other words, Theorem \ref{T2} generalizes the result obtained in \cite{KV} in the 2D situation.  In particular, the Gevrey bound \eqref{sc} and the radius of analyticity \eqref{radius} coincides with the ones derived in \cite{KV}. 
\end{rem}

\section{Global Well-posedness of the Nernst-Planck-Darcy system}\la{S4}

In this section, we prove Theorem \ref{T3}, concerning the global existence and uniqueness of solutions to  the NPD system: 

\beg{thm} \la{T3} (Global well-posedness and regularity of NPD) Let $T>0$ be arbitrary and $m \ge \fr{3}{2}$.  Assume that the initial ionic concentrations $c_i(0) \in H^m$ are nonnegative. Then the NPD system  described by \eqref{1}, \eqref{4} and \eqref{3} has a unique solution $(c_1, ..., c_n)$ on the time interval $[0,T]$ with the property that
\be 
(c_1, ..., c_n) \in (L^{\infty}(0,T; H^m))^n.
\ee
\end{thm}

 The proof of Theorem \ref{T3} is divided into several steps.

\textbf{Step 1. Existence of a unique local-in-time solution for $H^m$ initial datum.}
As in Section \ref{S2}, we focus on obtaining a priori estimates. We fix $m\ge \frac{3}{2}$. We apply $\l^m$ to \eqref{1}, multiply the resulting equation by $\l^m c_i$, and integrate by parts to obtain
\be
\begin{aligned}\la{dlmc}
\fr{1}{2}\fr{d}{dt}\|\l^m c_i\|_{L^2}+D_i\|\l^m\na c_i\|_{L^2}^2=&-(\l^m(u\cdot\na c_i),\l^m c_i)_{L^2}+D_iz_i(\l^m\na\cdot(c_i\na\Phi),\l^mc_i)_{L^2}\\
=&-([\l^m,u\cdot\na]c_i,\l^m c_i)_{L^2}-D_iz_i(\l^m(c_i\na\Phi),\l^m\na c_i)_{L^2}.
\end{aligned}
\ee
We estimate the first integral on the right hand side of \eqref{dlmc} as follows
\be
\begin{aligned}
\left|([\l^m,u\cdot\na]c_i,\l^mc_i)_{L^2}\right|\lesssim&(\|\l^m u\|_{L^2}\|\na c_i\|_{L^\infty}+\|{\na u}\|_{L^4}\|\l^m c_i\|_{L^4})\|\l^m c_i\|_{L^2}.
\end{aligned}
\ee
Then using $\|\na c_i\|_{L^\infty}\lesssim \|\l^m\na c_i\|_{L^2}$ for $m>1$ and $\|\l^m c_i\|_{L^4}\lesssim \|\l^m\na c_i\|_{L^2}$, we obtain
\be
\begin{aligned}\la{cc}
\left|([\l^m,u\cdot\na]c_i,\l^mc_i)_{L^2}\right|\le&\fr{D_i}{2}\|\l^m \na c_i\|_{L^2}^2+C(\|\l^m u\|_{L^2}^2+\|\na u\|_{L^4}^2)\|\l^m c_i\|_{L^2}^2.
\end{aligned}
\ee
Now, recalling from \eqref{3} that $u=-\mathbb{P}(\rho\na\Phi)$  where $\PP$ is the Leray projection onto  the space of divergence-free vectors, we have $\|\l^m u\|_{L^2}\le \|\l^m(\rho\na\Phi)\|_{L^2}$ and thus
\be
\begin{aligned}\la{cc1}
    \|\l^m u\|_{L^2}\lesssim  \|\l^m \rho\|_{L^2}\|\l^m\na\Phi\|_{L^2}\lesssim \sum_{j=1}^n\|\l^m c_j\|_{L^2}^2
\end{aligned}
\ee
which holds due to the Banach algebra property of $H^m$ for $m>1$. Similarly, we have
\be
\begin{aligned}\la{cc2}
\|\na u\|_{L^4} \lesssim \|\l^mu\|_{L^2}\lesssim \|\l^m\rho\|_{L^2}\|\l^m\na\Phi\|_{L^2}\lesssim \sum_{j=1}^n\|\l^m c_j\|_{L^2}^2
\end{aligned}
\ee
where the first inequality holds for $m\ge\frac{3}{2}$. Then, taking \eqref{cc1},\eqref{cc2} and returning to \eqref{cc}, we have
\be
\begin{aligned}\la{I1}
\left|([\l^m,u\cdot\na]c_i,\l^mc_i)_{L^2}\right|\le&\fr{D_i}{2}\|\l^m \na c_i\|_{L^2}^2+C\left(\sum_{j=1}^n\|\l^m c_j\|_{L^2}^2\right)^3.
\end{aligned}
\ee
Now we estimate the second integral on the right hand side of \eqref{dlmc}:
\be
\begin{aligned}\la{I2}
\left|D_iz_i(\l^m(c_i\na\Phi),\l^m\na c_i)_{L^2}\right|\le& \fr{D_i}{2}\|\l^m\na c_i\|_{L^2}^2+C\|\l^m(c_i\na\Phi)\|_{L^2}^2\\
\le&\fr{D_i}{2}\|\l^m\na c_i\|_{L^2}^2+C\|\l^m c_i\|_{L^2}^2\|\l^m\na\Phi\|_{L^2}^2\\
\le&\fr{D_i}{2}\|\l^m\na c_i\|_{L^2}^2+C\left(\sum_{j=1}^n\|\l^m c_j\|_{L^2}^2\right)^2.
\end{aligned}
\ee
Thus, summing \eqref{dlmc} in $i$ and using \eqref{I1} and \eqref{I2}, we obtain
\be
\fr{d}{dt}\left(\sum_{i=1}^n\|\l^m c_i\|_{L^2}^2\right)\lesssim \left(\sum_{i=1}^n\|\l^m c_i\|_{L^2}^2\right)^2+\left(\sum_{i=1}^n\|\l^m c_i\|_{L^2}^2\right)^3
\ee
from which it follows that for a short time $T^*$ depending on the $H^m$ norms of the initial data of $c_i$, there is a local solution of NPD on $[0,T^*)$ satisfying
\be\la{d'}
\sup_{t\in[0,T^*)}\sum_{i=1}^n\|\l^m c_i(t)\|_{L^2}^2\le 2\sum_{i=1}^n\|\l^m c_i(0)\|_{L^2}^2.
\ee
The uniqueness of local solutions follows from similar energy estimates. Suppose $(c_1, ..., c_n)$ and $(\bar c_1, ..., \bar c_n) $ are two local $H^m$ solutions to NPD on some common interval $[0,T*)$, where $T^*$ is taken small enough so that both solutions satisfy \eqref{d'}. Then, the differences $\tilde c_i=c_i-\bar c_i$ satisfy
\begin{align}
    \pa_t \tilde c_i+u\cdot \na\tilde c_i-D_i\D\tilde c_i=&D_iz_i\na\cdot(c_i\na\tilde\Phi+\tilde c_i\na\bar\Phi)-\tilde u\cdot\na\bar c_i.
\end{align}
Multiplying the above by $\tilde c_i$ and integrating by parts, we obtain exactly as in \eqref{tc}
\be
\begin{aligned}\la{tc'}
    \fr{1}{2}\fr{d}{dt}\|\tilde c_i\|_{L^2}^2+D_i\|\na\tilde c_i\|_{L^2}^2\le&\fr{D_i}{2}\|\na\tilde c_i\|_{L^2}^2+G_1(t)\left(\sum_{j=1}^n\|\tilde c_j\|_{L^2}^2+\|\tilde u\|_{L^2}^2\right)
\end{aligned}
\ee
with $G_1(t)$ uniformly bounded on $[0,T^*) $(c.f. \eqref{G1}). Next, we estimate, using $\tilde u=-\mathbb{P}(\rho\na\Phi-\bar\rho\na\bar\Phi)$
\be
\|\tilde u\|_{L^2}\lesssim \|\tilde\rho\|_{L^2}\|\na\Phi\|_{L^\infty}+\|\na\tilde\Phi\|_{L^2}\|\bar\rho\|_{L^\infty}\lesssim G_3(t)\sum_{j=1}^n \|\tilde c_j\|_{L^2}
\ee
where 
\be
G_3(t)=\|\na\Phi\|_{L^\infty}+\|\bar\rho\|_{L^\infty}
\ee
is uniformly bounded on $[0,T^*)$. Therefore, from \eqref{tc'}, we obtain after summing in $i$
\be
\fr{d}{dt}\sum_{i=1}^n\|\tilde c_i\|_{L^2}^2\lesssim G_1(t)\left(1+(G_3(t))^2\right)\sum_{i=1}^n\|\tilde c_i\|_{L^2}^2.
\ee
Then, applying Gronwall's inequality, we obtain
\be
\sup_{t\in[0,T^*)}\sum_{i=1}^n\|\tilde c_i(t)\|_{L^2}^2\le \sum_{j=1}^n\|\tilde c_i(0)\|_{L^2}^2e^{C\int_0^{T^*}G_1(s)(1+G_3(s))^2\,ds}.
\ee
Uniqueness follows from this last inequality.

\textbf{Step 2. Positivity of the ionic concentrations.} Suppose $(c_1, ..., c_n)$ is a solution of the NPD problem with the property that each ionic concentration $c_i$ belongs to the Lebesgue space $L^2(0,T; H^2)$. Then the positivity of $c_i(0)$ for $i \in \left\{1, ..., n\right\}$ is preserved for all positive times, that is $c_i(x,t) \ge 0$ for all $i \in \left\{1, ..., n\right\}$, for a.e. $x \in \TT^2$, and for all $t \in [0,T]$. The proof of this statement can be found in \cite{CI}. 

\textbf{Step 3.  A priori $L^2$ uniform-in-time bounds.} Suppose $(c_1, ..., c_n)$ is a solution of the NPD problem on the time interval $[0,T]$ such that $c_i(x,t) \ge 0$ for all $i \in \left\{1, ..., n\right\}$, for a.e. $x \in \TT^2$, and for all $t \in [0,T]$. Then there is a positive constant $\Gamma$ depending on the initial data, the time $T$, the parameters of the problem, and some universal constants, such that the following bound
\be \la{Gamma'}
\sup\limits_{0 \le t \le T}  \sum\limits_{i=1}^{n} \|c_i(t)\|_{L^2}  + \int_{0}^{T}\left(\|u\|_{L^2}^2 + \sum\limits_{i=1}^{n} D_i \|\na c_i(t)\|_{L^2}^2 \right) dt \le \Gamma
\ee 
holds. This can be obtained by following the proof of Proposition 2 in \cite{AI6}. We omit the details here. 

\textbf{Step 4. A priori $L^4$ uniform-in-time bounds.} We multiply the $i$-th ionic concentration equation \eqref{1} by $c_{i}^{3}$ and integrate the resulting equation over $\TT^2$. We obtain 
\be
\begin{split}
\frac{1}{4}\frac{d}{dt} \|c_i\|_{L^4}^4
+
\fr{3D_i}{4} \|\nabla c_{i}^2\|_{L^2}^2
&=
D_i z_i \int \nabla\cdot (c_i \nabla \Phi)c_{i}^{3}\,dx
=
-3D_i z_i \int c_i \nabla \Phi c_{i}^{2}\cdot \nabla c_i\,dx
\\&\lesssim \|c_{i}\|_{L^4}^{2} \|\nabla \Phi\|_{L^\infty} \|\nabla c_i^2\|_{L^2}
\\&\le
C\|c_{i}\|_{L^4}^{4} \|\nabla \Phi\|_{L^\infty}^2
+
\frac{3D_i}{4}\|\nabla c_i^2\|_{L^2}^2.
\end{split}
\ee Elliptic regularity yields the bound 
\be 
\|\nabla \Phi\|_{L^\infty} \lesssim \|\rho\|_{L^4} \lesssim \|\na \rho\|_{L^2}.
\ee Therefore, the $L^4$ norm of $c_i$ obeys
\be 
\fr{d}{dt} \|c_i\|_{L^4}^4 \lesssim \|c_i\|_{L^4}^4\|\na \rho\|_{L^2}^2.
\ee 
Integrating in time from $0$ to $t$ and using the boundedness of $\rho$ in $L^2(0,T;H^1)$ given by \eqref{Gamma'}, we conclude that the ionic concentrations are uniformly bounded on the space $L^{\infty}(0,T;L^4)$:
\be\la{c4}
\|c_i(t)\|_{L^4}\le \|c_i(0)\|_{L^4}\exp\left(C\int_0^t\|\na \rho(\tau)\|_{L^2}^2\,d\tau\right).
\ee
This finishes the proof of Step 4. 

\textbf{Step 5. A priori $H^1$ uniform-in-time bounds.} The gradient of the $i$-th ionic concentration evolves according to the energy equality 
\be \la{th3st5}
\fr{1}{2} \fr{d}{dt} \|\na c_i\|_{L^2}^2
+ D_i \|\Delta c_i\|_{L^2}^2
= (u \cdot \na c_i, \Delta c_i)_{L^2}
- D_i z_i (\na \cdot (c_i \na \Phi), \Delta c_i)_{L^2}.
\ee 
Using the divergence-free condition obeyed by the fluid velocity $u$ and elliptic regularity, we have
\be 
\|\na u\|_{L^2} 
\lesssim \|\na ^{\perp} \cdot u\|_{L^2}
\lesssim \|\na^{\perp} \rho \cdot \na \Phi\|_{L^2}
 \lesssim \|\na \rho\|_{L^2} \|\na \Phi\|_{L^{\infty}}
 \lesssim \|\na \rho\|_{L^2}\|\rho\|_{L^4},
\ee hence
\be 
\beg{aligned}
|(u \cdot \na c_i, \Delta c_i)_{L^2}|
& \lesssim\|u\|_{L^4} \|\na c_i\|_{L^4}\|\D c_i\|_{L^2}
\lesssim\|\rho\|_{L^4} \|\na \Phi\|_{L^\infty} \|\na c_i\|_{L^4}\|\Delta c_i\|_{L^2}\\
&\lesssim\|\rho\|_{L^4}\|\na\Phi\|_{L^\infty}\|\na c_i\|_{L^2}^{1/2}\|\D c_i\|_{L^2}^{3/2}\\
&\le \fr{D_i}{8} \|\Delta c_i\|_{L^2}^2 
+ C\|\rho\|_{L^4}^4\|\na\Phi\|_{L^\infty}^4 \|\na c_i\|_{L^2}^2.
\end{aligned}
\ee
Now we estimate
\be \la{th3st51}
\beg{aligned}
|D_i z_i (\na \cdot (c_i \na \Phi), \Delta c_i)_{L^2}|
& \lesssim \left(\|\na c_i\|_{L^2}\|\na \Phi\|_{L^{\infty}} + \|c_i\|_{L^4}\|\rho\|_{L^4} \right)\|\Delta c_i\|_{L^2}
\\& \lesssim \left(\|\na c_i\|_{L^2} \|\rho\|_{L^4} + \|c_i\|_{L^4} \|\rho\|_{L^4}\right) \|\Delta c_i\|_{L^2}
\\&\le \fr{D_i}{8} \|\Delta c_i\|_{L^2}^2 + C\|\rho\|_{L^4}^2 \|\na c_i\|_{L^2}^2 + C \|c_i\|_{L^4}^2 \|\rho\|_{L^4}^2
\end{aligned}
\ee 
Putting \eqref{th3st5}--\eqref{th3st51} together, we obtain the differential inequality
\be \la{nc}
\fr{d}{dt} \|\na c_i\|_{L^2}^2 + D_i \|\Delta c_i\|_{L^2}^2
\lesssim (\|\rho\|_{L^4}^4\|\na\Phi\|_{L^\infty}^4+\|\rho\|_{L^4}^2)\|\na c_i\|_{L^2}^2 +  \|c_i\|_{L^4}^2 \|\rho\|_{L^4}^2.
\ee Summing over $i \in \left\{1, ..., n \right\}$, and bounding the $L^4$ norm of the density $\rho$ by the sum of the $L^2$ norms of gradients of the ionic concentrations $c_i$, the latter energy inequality reduces to 
\be \la{bla}
\fr{d}{dt} \sum\limits_{i=1}^{n} \|\na c_i\|_{L^2}^2 + \sum\limits_{i=1}^{n} D_i \|\Delta c_i\|_{L^2}^2
\lesssim  \left(\|\rho\|_{L^4}^4\|\na\Phi\|_{L^\infty}^4+\|\rho\|_{L^4}^2 + \sum\limits_{i=1}^{n} \|c_i\|_{L^4}^2 \right)\sum\limits_{i=1}^{n} \|\na c_i\|_{L^2}^2 .
\ee 
By Gronwall's inequality, we obtain
\be
\begin{aligned}
\sum\limits_{i=1}^{n} \|\na c_i(t)\|_{L^2}^2 \le \sum\limits_{i=1}^{n} \|\na c_i(0)\|_{L^2}^2 \exp \left(\int_{0}^{t} \left(\|\rho(\tau)\|_{L^4}^4\|\na\Phi(\tau)\|_{L^\infty}^4+\|\rho(\tau)\|_{L^4}^2 + \sum\limits_{i=1}^{n} \|c_i(\tau)\|_{L^4}^2 \right) d\tau \right)
\end{aligned} 
\ee for all $t \ge 0$,
and so due to \eqref{c4} and \eqref{bla}, the ionic concentrations belong to the space
\be \la{lil2}
L^{\infty}(0,T; H^1) \cap L^2(0,T; H^2).
\ee
This completes the proof of Step 5.

\textbf{Step 6. A priori $H^m$ uniform-in-time bounds.} We multiply the equation obeyed by the $i$-th ionic concentration \eqref{1} by $\l^{2m} c_i$, and we integrate over the torus $\TT^2$. Integrating by parts, we obtain 
\be 
\fr{1}{2} \fr{d}{dt} \|\l^m c_i\|_{L^2}^2 
+ D_i \|\l^{m+1} c_i\|_{L^2}^2
=-(\l^{m-1} (u  \cdot\na c_i), \l^{m+1}  c_i)_{L^2}
- D_i z_i (\l^{m} (c_i \na \Phi), \l^{m} \na c_i)_{L^2}.
\ee Using fractional product estimates, the boundedness of the Leray projector on $L^4$, the Ladyzhenskaya's interpolation inequality, the Poincar\'e inequality, and the continuous embedding of $H^1$ in $L^{8}$, we estimate 
\be 
\beg{aligned}
\|\l^{m-1} u\|_{L^4} 
&= \|\l^{m-1} \PP (\rho \na \Phi)\|_{L^4}
\le \|\l^{m-1}(\rho \na \Phi)\|_{L^4}
\\&\lesssim \|\l^{m-1} \rho\|_{L^4}\|\na \Phi\|_{L^{\infty}}
+ \|\rho\|_{L^{8}}\|\l^{m-1} \na \Phi\|_{L^{8}}
\\&\lesssim \|\rho\|_{L^4} \|\l^{m} \rho\|_{L^2}
+  \|\na \rho\|_{L^2} \|\l^{m} \rho\|_{L^2}
\\&\lesssim \|\na \rho\|_{L^2} \|\l^{m} \rho\|_{L^2}.
\end{aligned}
\ee In view of the continuous embedding of $H^2$ in $L^{\infty}$ and the algebra property of $H^2$, we have
\be 
 \|u\|_{L^{\infty}}  \lesssim \|\Delta \PP (\rho \na \Phi)\|_{L^2} 
 \lesssim \|\Delta \rho\|_{L^2} \|\Delta \na \Phi\|_{L^2} \lesssim \|\Delta \rho\|_{L^2}\|\na\rho\|_{L^2}.
\ee 
Consequently, we bound the nonlinear term in $u$ as follows
\be 
\beg{aligned}
&|(\l^{m-1} (u \cdot \na c_i), \l^{m+1} c_i)_{L^2}|
\le \|\l^{m+1} c_i\|_{L^2} \|\l^{m-1} (u \cdot \na c_i)\|_{L^2}
\\&\lesssim \|\l^{m+1} c_i\|_{L^2} \left(\|\l^{m-1} u\|_{L^4} \|\na c_i\|_{L^4} + \|u\|_{L^{\infty}} \|\l^{m} c_i\|_{L^2} \right)
\\&\lesssim \|\l^{m+1} c_i\|_{L^2} \left(\|\Delta c_i\|_{L^2} \|\na \rho\|_{L^2} \|\l^{m} \rho\|_{L^2}  +\|\Delta \rho\|_{L^2} \|\na \rho\|_{L^2} \|\l^{m} c_i\|_{L^2}  \right).
\end{aligned}
\ee Now we estimate 
\be 
\beg{aligned}
&|(\l^{m} (c_i \na \Phi), \l^{m} \na c_i)_{L^2}|
\lesssim\|\l^{m+1} c_i\|_{L^2} \left(\|\l^m c_i\|_{L^2} \|\na \Phi\|_{L^{\infty}} + \|c_i\|_{L^4}\|\l^{m} \na \Phi\|_{L^4} \right)
\\&\lesssim \|\l^{m+1} c_i\|_{L^2} \left(\|\l^m c_i\|_{L^2}  \|\na \rho\|_{L^2} + \|c_i\|_{L^4}\|\l^{m} \rho\|_{L^2} \right).
\end{aligned}
\ee This gives the energy inequality
\be 
\beg{aligned}
&\fr{d}{dt} \|\l^{m} c_i\|_{L^2}^2 
+ D_i \|\l^{m+1} c_i\|_{L^2}^2
\\&\lesssim\|\Delta c_i\|_{L^2}^2  \|\na \rho\|_{L^2}^2 \|\l^{m} \rho\|_{L^2}^2  + \|\Delta \rho\|_{L^2}^2 \|\na \rho\|_{L^2}^2 \|\l^{m} c_i\|_{L^2}^2 +  \|\na \rho\|_{L^2}^2  \|\l^m c_i\|_{L^2}^2 + \|c_i\|_{L^4}^2\|\l^{m} \rho\|_{L^2.}^2
\end{aligned}
\ee Finally, we use the triangle inequality to bound
\be 
\|\l^{m} \rho\|_{L^2}^2  \lesssim \sum\limits_{j=1}^{n} \|\l^{m} c_j\|_{L^2}^2,
\ee sum over all indices $i \in \left\{1, ..., n \right\}$, integrate in time from $0$ to $t$,  exploit the regularity obtained in \eqref{lil2}, and conclude that the ionic concentrations satisfy
\be\la{cBt}
\sum_{i=1}^n\|\l^m c_i(t)\|_{L^2}^2 \,d\tau\le e^{L(t)}\sum_{i=1}^n\|\l^m c_i(0)\|_{L^2}^2
\ee
where
\be\la{B(t)}
L(t)=C\left(\sup_{\tau\in[0,t)}  \|\na \rho(\tau)\|_{L^2}^2 \right)\left(\sum_{i=1}^n\int_0^t  \|\Delta c_i(\tau)\|_{L^2}^2\,d\tau\right)+\int_0^t \|\na \rho(\tau)\|_{L^2}^2 \,d\tau+\sum_{i=1}^n\int_0^t\|c_i(\tau)\|_{L^4}^2\,d\tau.
\ee
Thus for all $i \in \left\{1, ..., n\right\}$, $c_i$ belongs to the space $L^{\infty}(0, T; H^m)$.

\textbf{Step 7. Extension of the local solution.} The local solution can be extended to the time interval $[0,T]$, a fact that follows from the uniform-in-time boundedness in $H^m$ obtained in Step 6. This ends the proof of Theorem~\ref{T3}.

\section{Gevrey Regularity of the Nernst-Planck-Darcy system} \la{S5}

In this section, we address the propagation of the Gevrey regularity for the NPD system:

\beg{thm} \la{T4} (Global analyticity of NPD) Let $T>0$ be arbitrary and $m > 2$. Assume that the initial ionic concentrations $c_i(0) \in H^m$ are nonnegative. Then the NPD system  described by \eqref{1}, \eqref{4} and \eqref{3} has a unique solution $(c_1, ..., c_n)$ on [0,T] such that for each time $t \in (0,T)$, the functions $c_1, ..., c_n$ are real analytic in the spatial variable with uniform radius of analyticity 
\be 
\tau(t) = \fr{1}{2} \min\left\{D_1, ..., D_n \right\} \min\{t,T_0/2\}
\ee
where $T_0>0$ depends only on parameters, the $H^m$ norms of $c_i(0)$, and lower regularity Sobolev norms of $c_i$ up to time $T$.
\end{thm}

In order to prove Theorem \ref{T4}, we take the scalar products in $\mathcal{D}(e^{\tau \Lambda})$ of the equation obeyed by the ionic concentration $c_i$ with $\l^{2m}c_i$.  We obtain
\be \la{energyD}
\beg{aligned} 
&\fr{1}{2} \fr{d}{dt} \left[\sum\limits_{i=1}^{n} \|e^{\tau \l} \l^m c_i\|_{L^2}^2 \right] 
+ \sum\limits_{i=1}^{n} D_i \|e^{\tau \l} \l^{m+1} c_i\|_{L^2}^2 
-\tau'(t) \sum\limits_{i=1}^{n} \|e^{\tau \l} \l^{m + \fr{1}{2}} c_i\|_{L^2}^2  
\\&\quad=- \sum\limits_{i=1}^{n} (e^{\tau \l} \l^m (u \cdot \na (c_i - \bar{c_i}) ), e^{\tau \l} \l^m(c_i - \bar{c}_i))_{L^2}
+ \sum\limits_{i=1}^{n} D_i z_i (e^{\tau \l} \l^{m} \na \cdot [(c_i - \bar{c}_i) \na \Phi], e^{\tau \l} \l^m (c_i - \bar{c}_i))_{L^2} 
\\&\quad\quad\quad\quad+ \sum\limits_{i=1}^{n} D_i z_i \bar{c}_i (e^{\tau \l}\l^m \Delta \Phi, e^{\tau \l} \l^m (c_i - \bar{c}_i))_{L^2}. 
\end{aligned} \ee

In order to control the nonlinear terms, we need the following lemmas.

\beg{lem} \la{L6} Let $i \in \left\{1, ..., n \right\}$. The following estimate
\be   \la{nonlin6}
\beg{aligned}
& \sum\limits_{i=1}^{n} |(e^{\tau \l} \l^{m} (u \cdot \na (c_i - \bar{c}_i )), e^{\tau \l} \l^m (c_i - \bar{c}_i ))_{L^2}| 
\\&\quad\quad\leq \sum\limits_{i=1}^{n} \fr{D_i}{8} \|e^{\tau \l} \l^{m+1}  c_i\|_{L^2}^2
+ C \left(\sum\limits_{i=1}^{n} \|e^{\tau \l} \l^m c_i\|_{L^2}^2 \right)^3
\end{aligned} \ee holds for any $m >2 $.
\end{lem}

\textbf{Proof.} We set $c_i$ and $c_i^*$ as in \eqref{cfourier} and \eqref{cstarfourier} respectively. We write the Fourier expansion of $\rho$
\be 
\rho = \sum\limits_{v \in \ZZ^2 \setminus \left\{0\right\}} \rho_v e^{i v \cdot x},
\ee from which we obtain the Fourier expansion of $\rho \na \Phi$ 
\be 
\rho \na \Phi = i \sum\limits_{v, z \in \ZZ^2 \setminus \left\{0\right\}} \fr{z}{|z|^2} \rho_z \rho_v e^{i(v+z) \cdot x}
= i \sum\limits_{j \in \ZZ^2 \setminus \left\{0\right\}} \left[\sum\limits_{z \in \ZZ^2 \setminus \left\{0, j\right\}} \fr{z}{|z|^2} \rho_{z} \rho_{j-z} \right]e^{ij \cdot x}.
\ee
In view of the divergence free-condition obeyed by $u$, we have 
\be 
u = - \PP(\rho \na \Phi)
\ee where $\PP$ is the Leray projection onto the space of divergence-free vectors. Consequently, the Fourier series of $u$ is given by 
\be 
u = - \sum\limits_{j \in \ZZ^2 \setminus \left\{0\right\}} \left[v_j - (v_j \cdot j)\fr{j}{|j|^2} \right] e^{ij \cdot x}
\ee where
\be 
v_j = \sum\limits_{z \in \ZZ^2 \setminus \left\{0,j\right\}} i\fr{z}{|z|^2} \rho_z \rho_{j-z}.
\ee For $j \in \ZZ^2 \setminus \left\{0\right\}$, we denote by $u_j$ the Fourier coefficients of $u$, that is 
\be 
u_j = - \left[v_j - (v_j \cdot j)\fr{j}{|j|^2} \right],
\ee and we obtain 
\be 
\beg{aligned}
&(e^{\tau \l} \l^{m} (u \cdot \na (c_i - \bar{c}_i )), e^{\tau \l} \l^m (c_i - \bar{c}_i ))_{L^2}
= 4\pi^2 i \sum\limits_{j+k+l=0} (u_j \cdot k) (c_i)_k (c_i)_l |l|^{2m} e^{2\tau l} 
\\&\quad\quad= 4\pi^2 \sum\limits_{j+k+l=0} \sum\limits_{z \in \ZZ^2 \setminus \left\{0,j\right\}} \left[\fr{z \cdot k}{|z|^2} - \fr{z \cdot j}{|z|^2} \fr{j \cdot k}{|j|^2} \right] \rho_z \rho_{j-z} (c_i)_k (c_i)_l |l|^{2m} e^{2\tau l} 
\\&\quad\quad= 4\pi^2 \sum\limits_{z + \tilde{z}+k+l=0} \left[\fr{z \cdot k}{|z|^2} - \fr{z \cdot (z+\tilde{z})}{|z|^2} \fr{(z+\tilde{z}) \cdot k}{|z + \tilde{z}|^2} \right] \rho_z \rho_{\tilde{z}} (c_i)_k (c_i)_l |l|^{2m} e^{2\tau l}.
\end{aligned} \ee We point out that the last sum is over all indices $z, \tilde{z}, k, l \in \ZZ^2$ such that $z + \tilde{z} + k + l = 0, \tilde{z} + z \ne 0, z \ne 0, k \ne 0$, $l \ne 0$ and $\tilde{z} \ne 0$. Using the fact that $z + \tilde{z}+k+l=0$, we have $|l| \le |k| + |z| + |\tilde{z}|$, hence 
\be \la{lem61}
\beg{aligned} 
&|(e^{\tau \l} \l^{m} (u \cdot \na (c_i - \bar{c}_i )), e^{\tau \l} \l^m (c_i - \bar{c}_i ))_{L^2}| 
\\&\quad\quad\lesssim  \sum\limits_{z + \tilde{z}+k+l=0} |\rho_z^*| |\rho_{\tilde{z}}^*| |(c_i)_k^*| |(c_i)_l^*| |k| |z|^{-1} |l|^{m+1} (|k|^{m-1} + |z|^{m-1} + |\tilde{z}|^{m-1})
\end{aligned} \ee where $\rho_{z}^* = e^{\tau|z|} \rho_z$ and $\rho_{\tilde{z}}^* = e^{\tau |\tilde{z}|} \rho_{\tilde{z}}$. We estimate the following three terms
\be 
\beg{aligned}
&\sum\limits_{z + \tilde{z}+k+l=0} |\rho_z^*| |\rho_{\tilde{z}}^*| |(c_i)_k^*| |(c_i)_l^*| |k|^m |z|^{-1} |l|^{m+1} \lesssim \||l|^{m+1} (c_i)_{l}^*\|_{\ell^2} \||k|^m (c_i)_{k}^*\|_{\ell^2} \||z|^{-1} \rho_{z}^*\|_{\ell^1} \|\rho_{\tilde{z}}^*\|_{\ell^{\infty}} 
\\&\quad\quad\quad\quad\lesssim \||l|^{m+1} (c_i)_{l}^*\|_{\ell^2} \||k|^m (c_i)_{k}^*\|_{\ell^2} \||z|^{-1} \rho_{z}^*\|_{\ell^1} \|\rho_{\tilde{z}}^*\|_{\ell^{1}} 
\\&\quad\quad\quad\quad\lesssim \|\l^{m+1} e^{\tau \l} c_i\|_{L^2} \|\l^{m} e^{\tau \l} c_i\|_{L^2}\|\l^{\epsilon} e^{\tau \l} \rho \|_{L^2} \|\l^{1+ \epsilon} e^{\tau \l} \rho \|_{L^2},
\end{aligned} \ee
\be 
\beg{aligned}
&\sum\limits_{z + \tilde{z}+k+l=0} |\rho_z^*| |\rho_{\tilde{z}}^*| |(c_i)_k^*| |(c_i)_l^*| |k| |z|^{m-2} |l|^{m+1} 
 \lesssim \||l|^{m+1} (c_i)_{l}^*\|_{\ell^2} \||z|^{m-2} \rho_z^*\|_{\ell^2} \||k|(c_i)_k^*\|_{\ell^1} \|\rho_{\tilde{z}}^*\|_{\ell^{\infty}} 
\\&\quad\quad\quad\quad\lesssim \|\l^{m+1} e^{\tau \l} c_i\|_{L^2} \|\l^{m-2} e^{\tau \l} \rho\|_{L^2}\|\l^{2+ \epsilon} e^{\tau \l} c_i \|_{L^2} \|\l^{1+ \epsilon} e^{\tau \l} \rho \|_{L^2},
\end{aligned} \ee
\be \la{lem62}
\beg{aligned} 
&\sum\limits_{z + \tilde{z}+k+l=0} |\rho_z^*| |\rho_{\tilde{z}}^*| |(c_i)_k^*| |(c_i)_l^*| |k| |z|^{-1} |l|^{m+1} |\tilde{z}|^{m-1} 
\lesssim \||l|^{m+1} (c_i)_{l}^*\|_{\ell^2} \||\tilde{z}|^{m-1} \rho_{\tilde{z}}^*\|_{\ell^2} \||z|^{-1} \rho_{z}^*\|_{\ell^{\infty}} \||k|(c_i)_{k}^*\|_{\ell^{1}} 
\\&\quad\quad\quad\quad\lesssim \|\l^{m+1} e^{\tau \l} c_i\|_{L^2} \|\l^{m-1} e^{\tau \l} \rho\|_{L^2}\|\l^{2+ \epsilon} e^{\tau \l} c_i \|_{L^2} \|\l^{\epsilon} e^{\tau \l} \rho \|_{L^2}
\end{aligned} \ee
for any $m > 0$ and any $\epsilon > 0$. Putting \eqref{lem61}--\eqref{lem62} together and summing over all $i \in \left\{1, ..., n\right\}$, we obtain \eqref{nonlin6} for any $m > 2$.

\beg{lem} \la{L7}  Let $i \in \left\{1, ..., n \right\}$.  The following estimate 
\be \la{nonlin7}
\beg{aligned} 
&\sum\limits_{i=1}^{n} D_i |z_i| |(e^{\tau \l} \l^{m} \na \cdot [(c_i - \bar{c}_i) \na \Phi], e^{\tau \l} \l^m (c_i - \bar{c}_i))_{L^2}| 
\\&\quad\quad\le \sum\limits_{i=1}^{n} \fr{D_i}{8} \|e^{\tau \l} \l^{m+1}  c_i\|_{L^2}^2
+ C \left(\sum\limits_{i=1}^{n} \|e^{\tau \l} \l^m c_i\|_{L^2}^2 \right)^2
\end{aligned} \ee holds for any $m > 1.$ 
\end{lem}

\textbf{Proof.} Under the same setting introduced in the proof of Lemma \ref{L3}, we have 
\be 
\beg{aligned} 
&\left|(e^{\tau \l} \l^{m} \na \cdot [(c_i - \bar{c}_i) \na \Phi], e^{\tau \l} \l^m (c_i - \bar{c}_i))_{L^2}\right|
= 4\pi^2 \left|\sum\limits_{j + k + l = 0} (c_i)_{j} \rho_k (c_i)_l \fr{l \cdot k}{|k|^2}|l|^{2m} e^{2\tau |l|}\right| 
\\&\quad\quad\lesssim \sum\limits_{j+k+l = 0} |(c_i)_j^*| |\rho_k^*| |(c_i)_l^*| |l|^{m+1} |k|^{-1} (|k|^m + |j|^m) 
\\&\quad\quad\lesssim \||l|^{m+1}(c_i)_{l}^*\|_{\ell^2} \left[\||k|^{m-1}\rho_k^*\|_{\ell^2} \|(c_i)_j^*\|_{\ell^1} + \||j|^m (c_i)_j^*\|_{\ell^2} \||k|^{-1}\rho_k^*\|_{\ell^1} \right],
\end{aligned} \ee which gives the desired estimate \eqref{nonlin7} after applications of Plancherel's identity and Young's inequality. 

\beg{lem} \la{L8} Let $i \in \left\{1, ..., n \right\}$. The following estimate
\be \la{nonlin8}
\beg{aligned} 
\sum\limits_{i=1}^{n} D_i z_i \bar{c}_i |(e^{\tau \l}\l^m \Delta \Phi, e^{\tau \l} \l^m (c_i - \bar{c}_i))_{L^2}| 
\le C\sum\limits_{i=1}^{n} \|e^{\tau \l} \l^m c_i\|_{L^2}^2
\end{aligned} \ee holds for any $m > 0 $. Here $C$ is a positive constant depending on the parameters of the problem and the initial spatial average of the ionic concentrations. 
\end{lem}

\textbf{Proof.} The proof follows from a direct application of the Cauchy-Schwarz inequality.

Now we complete the proof of Theorem \ref{T4}. For fixed $T>0$, it follows from \eqref{cBt} that
\be\la{BT}
\sup_{t\in[0,T]}\sum_{i=1}^n\|\l^m c_i(t)\|_{L^2}^2 \,d\tau\le e^{L(T)}\sum_{i=1}^n\|\l^m c_i(0)\|_{L^2}^2=:L_T.
\ee
Now, choosing $\tau(t) = \fr{1}{2} \min\left\{D_1, ..., D_n \right\} t$, and putting \eqref{energyD}, \eqref{nonlin6}, \eqref{nonlin7} and \eqref{nonlin8} together, we obtain the energy inequality
\be \la{33}
\fr{d}{dt} \left[\sum\limits_{i=1}^{n} \|e^{\tau \l} \l^m c_i\|_{L^2}^2 \right]\lesssim \left[1+ \sum\limits_{i=1}^{n} \|e^{\tau \l} \l^m c_i\|_{L^2}^2 \right]^3 
\ee from which we deduce the boundedness of the Gevrey norm   
\be 
\sum\limits_{i=1}^{n} \|e^{\tau \l} \l^m c_i\|_{L^2}^2
\le2\left(1+ \sum\limits_{i=1}^{n} \|e^{\tau(0)\l}\l^m c_i(0)\|_{L^2}^2\right)= 2\left(1+ \sum\limits_{i=1}^{n} \|\l^m c_i(0)\|_{L^2}^2\right)
\ee on a short time interval $[0,T_0]$, where $T_0$ depends on the initial magnitude in $L^2$ of $e^{\tau\l}\l^mc_i$, which is simply the $H^m$ norm of $c_i$. In particular, we may choose $T_0$ to just be a function of parameters and $L_T$ from \eqref{BT}. This yields local Gevrey regularity.
To extend Gevrey regularity to a longer interval, we take this time $\tau(t)=\fr{1}{2}\min\{D_1,...,D_n\}(t-\fr{T_0}{2})$ and obtain from \eqref{33} the following bound
\be 
\sum\limits_{i=1}^{n} \|e^{\tau \l} \l^m c_i\|_{L^2}^2
\le2\left(1+ \sum\limits_{i=1}^{n} \|e^{\tau(T_0/2)\l}\l^m c_i(T_0/2)\|_{L^2}^2\right)= 2\left(1+ \sum\limits_{i=1}^{n} \|\l^m c_i(T_0/2)\|_{L^2}^2\right)
\ee
on the interval $[T_0/2,T_0/2+T_1]$, where $T_1$ depends on the $H^m$ norms of $c_i$ at time $T_0/2$. However, due to the uniform bound \eqref{BT}, we may choose $T_1=T_0$. Thus, we have Gevrey regularity on the interval $[0,3T_0/2]$, and in the same manner we extend to the whole interval $[0,T]$. Based on this construction, we see that the radius of analyticity $\tau$ is bounded below by $\fr{1}{2}\min\{D_1,...,D_n\}t$ on the interval $[0,T_0]$ and then bounded below by $\fr{1}{2}\min\{D_1,...,D_n\}\fr{T_0}{2}$ on $[T_0,T]$.
This ends the proof of Theorem \ref{T4}.

\section{Acknowledgments} 
WW was partially supported by an AMS-Simons travel grant.

\end{document}